\documentclass[amscd,amssymb,verbatim,10pt,letterpaper]{rmmcart}
\usepackage{latexsym, amssymb, graphics, setspace,amscd}
\usepackage{euscript}
\usepackage[small,nohug,,heads=LaTeX]{diagrams}
\diagramstyle[labelstyle=\scriptstyle]

    \newcommand{\ie}{{\em i.e.}}
    \newcommand{\eg}{{\em e.g.}}

    \newtheorem{thm}{Theorem}[section]
    \newtheorem{prop}[thm]{Proposition}
    \newtheorem{lem}[thm]{Lemma}
    \newtheorem{cor}[thm]{Corollary}

    \theoremstyle{definition}

    \newtheorem{conj}[thm]{Conjecture}

    \theoremstyle{remark}

    \newtheorem{rem}[thm]{Remark}

    \newcommand{\numbtext}[4]{\begin{equation}  \label{#1}
				  \left#2
				      \begin{array}{cc}
					  {} \parbox{4in}{#4} & {}
				      \end{array}
				  \right#3
			      \end{equation}}

    \newcommand{\fn}[3]{#1 \colon #2 \rightarrow #3}

    \newcommand{\fna}[2]{#1 \rightarrow #2}

    \newcommand{\fnai}[2]{#1 \hookrightarrow #2}

    \newcommand{\fnl}[3]{#2 \stackrel{#1}{\rightarrow} #3}

    \newcommand{\Spec}{\operatorname{Spec}}

     \newcommand{\Hilb}{\operatorname{Hilb}}

    \newcommand{\Hom}[3]{\operatorname{Hom}_{#1}(#2,#3)}

    \newcommand{\gf}{K}   

    \newcommand{\Span}{\operatorname{Span}}

\begin{document}

\title[Some Elementary Components]{Some Elementary Components of the Hilbert Scheme of Points}

\author[M. Huibregtse]{Mark E. Huibregtse}

\address{Department of Mathematics and Computer Science\\
         Skidmore College\\
         Saratoga Springs, New York 12866}

\email{mhuibreg@skidmore.edu}


\date{\today} 


\newcommand{\leadmon}{\operatorname{LM}} 
\newcommand{\trailmon}{\operatorname{TM}} 
\newcommand{\tarmons}{\operatorname{T}} 
\newcommand{\tansp}{\EuScript{T}}
\newcommand{\specl}{\mathcal{S}}
\newcommand{\SLI}{distinguished}  
\newcommand{\SLICap}{Distinguished}
\newcommand{\eff}{efficient}


\begin{abstract}
    
Let $K$ be an algebraically closed field of characteristic $0$, and let $H^{\mu}_{\mathbb{A}^n_{\gf}}$ denote the Hilbert scheme of $\mu$ points of $\mathbb{A}^n_{K}$.  An \textbf{elementary component} $E$ of $H^{\mu}_{\mathbb{A}^n_{\gf}}$ is an irreducible component such that every $K$-point $[I]$ $\in$ $E$ represents a length-$\mu$ closed subscheme $\Spec(\gf[x_1,\dots,x_n]/I)$ $\subseteq$ $\mathbb{A}^n_{K}$ that is supported at one point.  Iarrobino and Emsalem gave the first explicit examples (with $\mu > 1$) of elementary components in \cite{Iarrob-Emsalem}; in their examples, the ideals $I$ were homogeneous (up to a change of coordinates corresponding to a translation of $\mathbb{A}^n_{K}$).  We generalize their construction to obtain new examples of elementary components.
\end{abstract}

\maketitle

\section{Introduction} \label{sec:Intro}


Let $\gf$ be an algebraically closed field of characteristic $0$\footnote{This hypothesis is used explicitly in Proposition \ref{prop:GeneralZProp} and subsequent results that depend thereon, and implicitly in the computer computations, which are (with one exception) done in characteristic $0$.}, and let $R$ denote the polynomial ring $\gf[x_1, \dots, x_n]$ $=$ $\gf[\mathbf{x}]$, with $n$ $\geq$ $3$.  The first explicit examples of finite algebras with ``small tangent space'' (or ``generic'' algebras) of $\gf$-dimension $\mu > 1$ were given by Iarrobino and Emsalem in \cite{Iarrob-Emsalem}.  These algebras have the form $A$ $=$ $R/I$, where $I$ $\subseteq$ $R$ is an ideal of finite colength $\mu$ that is generated by a list of sufficiently general homogeneous polynomials  $g_j$, $1 \leq j \leq \lambda$, of degree $r$, and so vanishes at a single point (the origin) of $\mathbb{A}^n_{\gf}$.  
The point $[I]$ $\in$ $\Hilb^{\mu}_{\mathbb{A}^n_{\gf}}$ corresponding to the ideal $I$ has a small tangent space in the sense that the tangent directions at $[I]$ correspond to deformations of $I$ to ideals $I'$ of the same ``type,'' obtained either by varying the coefficients of the generators $g_j$ or by translating the subscheme $\Spec(R/I)$ in $\mathbb{A}^n_{\gf}$; in particular, all of the $I'$ vanish at a single point.   Accordingly, the point $[I]$ is a simple point on an \textbf{elementary component} of the Hilbert scheme; that is, a component such that every point on it parameterizes a subscheme concentrated at a single point \cite[p.\ 148]{Iarrob:Sitges1983}.  If $\mu$ $>$ $1$, it is clear that an elementary component must be different from the \textbf{principal component}, which contains the points corresponding to reduced subschemes of length $\mu$.  (Note that if $\mu > 1$ for an elementary component, then $\mu \geq 8$; see \cite[Th.\ 1.1]{CartwrightErmanVelascoViray:HilbSchOf8Pts}.)  The purpose of this paper is to generalize the construction in \cite{Iarrob-Emsalem} to produce new examples of generic algebras $R/I$ (or elementary components of $\Hilb^{\mu}_{\mathbb{A}^n_{\gf}}$). 

\begin{rem} \label{rem:History}
    Since a $0$-dimensional closed subscheme of $\mathbb{A}^n_{\gf}$ can be written as a disjoint union of subschemes supported at single points of $\mathbb{A}^n_{\gf}$, one sees that elementary components are the ``building blocks'' of irreducible components of $H^{\mu}_{\mathbb{A}^n_{\gf}}$.  Hence, Iarrobino's demonstration that $H^n_{\mathbb{A}^{\mu}_{\gf}}$ is in general reducible \cite{Iarrob:ReducibilityOfHilb} already implies the existence of non-trivial elementary components (see \cite{Iarrob:NumberOfGenericSingularities}.)  As previously noted, the first explicit examples were given by Iarrobino and Emsalem in \cite{Iarrob-Emsalem}.  Employing a different approach, Shafarevich gave further examples in \cite{Shaf:DefsOfCommAlgebrasOfClass2}.
\end{rem} 

The present paper is but one small contribution to the voluminous, diverse, and rapidly increasing literature on components of Hilbert schemes of points; for example, see  \cite{BorgesDosSantosHenniJardim:CommutingMatrices}, 
\cite{CartwrightErmanVelascoViray:HilbSchOf8Pts}, \cite{CasnatiJelisiejewNotari:GorensteinLociRayFamilies}, \cite{ErmanAndVelasco:SyzygeticApproach}, and the references contained therein.

\subsection{Iarrobino-Emsalem example} \label{ssec:I-Eeg} To set the stage for our generalization, we describe more fully Iarrobino and Emsalem's first example using our notation and terminology.  The ideal $I$ $\subseteq$ $R = \gf[x_1, x_2, x_3, x_4]$ is generated by quadratic forms 
\[
    g_j = m_j - N_j,\ 1 \leq j \leq 7,
\]
where $m_j$ denotes the $j$-th monomial in the list of ``leading'' monomials
\[
    x_1^2,\, x_1 x_2,\, x_1 x_3,\, x_1 x_4,\, x_2^2,\, x_2 x_3,\, x_2 x_4, 
\]
 and 
\[
    N_j = \sum_{i=0}^{2}(c_{i}\cdot x_3^{i}x_4^{2-i})
\]
is a $\gf$-linear combination of the ``trailing'' monomials of degree $2$ in the ``back variables'' $x_3,\, x_4$.  When the coefficients $c_{i}$ are sufficiently general, one can show that all the monomials of degree $3$ belong to $I$; consequently, $I$ has finite colength with zero-set concentrated at the origin, and one sees easily that the order ideal 
\[
    \mathcal{O} = \{1,\, x_1,\, x_2,\, x_3,\, x_4,\, x_3^2,\, x_3 x_4,\, x_4^2 \}
\]
 is a $\gf$-basis of the quotient $R/I$.  Therefore, $[I]$ $\in$ $\Hilb^8_{\mathbb{A}^4_{\gf}}$; moreover, $I$ is in the ``border basis scheme'' $\mathbb{B}_{\mathcal{O}}$ $\subseteq$ $\Hilb^8_{\mathbb{A}^4_{\gf}}$ (see \cite[Secs.\ 2, 3]{KreutzerAndRobbiano1:DefsOfBorderBases}); we recall briefly the basics of border basis schemes in Section \ref{sec:borderBases}.  

The ideal $I$ can be ``deformed'' in two ways: the $7 \cdot 3$ $=$ $21$ coefficients defining the $g_j$ can be tweaked, and the  ideal (or corresponding subscheme) can be translated in four independent directions in $\mathbb{A}^4_{\gf}$; this shows that the point $[I]$ lies on a locus of dimension at least 25 consisting of points $[I']$ such that $I'$ is supported at one point.  On the other hand, the dimension of the tangent space $\tansp_{[I]}$ $=$ $\Hom{R}{I}{R/I}$ can be computed, and one finds that this dimension is 25.  From this it follows that $[I]$ is a smooth point on an elementary component of dimension 25 in $\Hilb^{8}_{\mathbb{A}^4_{\gf}}$.  In general, we say that an ideal $I$  such that $[I]$ is a smooth point on an elementary component of $\Hilb^{\mu}_{\mathbb{A}^n_{\gf}}$ is \textbf{generic}.

\begin{rem} \label{rem:HilbFuncAndShaf}
    The Hilbert function of the subschemes $\operatorname{Spec}(R/I)$ in the example just discussed is $(1,4,3,0)$.  Shafarevich's results in \cite{Shaf:DefsOfCommAlgebrasOfClass2} imply that the analogously-constructed ideals corresponding to the Hilbert function $(1,5,3,0)$ are also generic.  In the next section we describe our generalization of Iarrobino and Emsalem's construction, which yields generic ideals corresponding to the Hilbert function $(1,5,3,4,0)$.
\end{rem}

\subsection{A generalization} \label{ssec:AGen}

We now describe our ``smallest'' example of a generic ideal $I$; it is very similar to the Iarrobino-Emsalem example just discussed, except that the leading and trailing monomials have different degrees and the embedding dimension is $5$. (A more complete description, including a link to a \textit{Mathematica} \cite{Mathematica} notebook containing the computational details, is given in Section \ref{ssec:(1,5,3,4)}.)  The leading monomials are the $12$ monomials of degree $2$ in $R$ $=$ $\gf[x_1,\dots,x_5]$ that involve at least one of the ``front variables'' $x_1, x_2, x_3$:
\[
  \leadmon = \left\{ \begin{array}{c}
      x_1^2,\, x_1 x_2,\, x_1 x_3,\, x_1 x_4,\, x_1 x_5,\, x_2^2,\, x_2 x_3,\\ x_2 x_4,\, x_2 x_5, x_3^2,\, x_3 x_4,\, x_3 x_5
  \end{array} \right\},
\]
and the trailing monomials are the four monomials of degree $3$ in the ``back variables'' $x_4, x_5$:
\[
    \trailmon = \{x_4^3,\, x_4^2 x_5,\, x_4 x_5^2,\, x_5^3\}.
\] 
The ideal is generated by polynomials
\[
    g_j = m_j - N_j,\ 1 \leq j \leq 12,
\]
where $m_j$ is the $j$-th leading monomial and 
\[
    N_j = \sum_{i=0}^{3} c_{i}x_4^{i}x_5^{3-i} \in \Span_{\gf}(\trailmon)
\]
is a form of degree $3$ in $x_4, x_5$.  If the $g_j$ are sufficiently general, it can be shown that every monomial of degree $4$ is in $I$, and that the quotient $R/I$ has $\gf$-basis the order ideal
\[
    \mathcal{O} = \left\{1,\, x_1,\, x_2,\, x_3,\, x_4,\, x_5,\, x_4^2,\, x_4x_5,\, x_5^2,\, x_4^3,\, x_4^2x_5,\, x_4x_5^2,\, x_5^3\right\},
\]
so 
\[
    [I] \in \mathbb{B}_{\mathcal{O}} \subseteq \Hilb^{13}_{\mathbb{A}^5_{\gf}}
\]
and the Hilbert function is $(1,5,3,4,0)$.
 
As shown in general in Sections \ref{sec:idealFam} and \ref{sec:derivMap}, there are (at least) three ways that the ideal $I$ can be deformed without changing its ``type,'' or the fact that its zero-set consists of one point, and that these give independent tangent directions at $[I]$: 
\begin{itemize}
  \item The $4 \cdot 12 = 48$ coefficients $c_{ij}$ can be tweaked;
  \item The ideal can be translated in $\mathbb{A}^5_{\gf}$;
  \item The ideal can be pulled back via automorphisms of $\mathbb{A}^5_{\gf}$ defined by coordinate changes of the form $x_{\alpha}$ $\mapsto$ $x_{\alpha} + c_{\alpha,\beta}\cdot x_{\beta}$, $x_{\beta}$  $\mapsto$ $x_{\beta}$, where $1 \leq \alpha \leq 3$, $4 \leq \beta \leq 5$, $c_{\alpha,\beta} \in \gf$.
\end{itemize} 
(Translation also involves pulling back the ideal via an automorphism of $\mathbb{A}^5_{\gf}$, and so the second and third deformation methods are treated in a uniform way in the body of the paper.)  Therefore, $[I]$ lies on a locus of dimension at least $48 + 5 + 3 \cdot 2$ $=$ $59$ consisting of points $[I']$ such that the ideal $I'$ is supported at one point.  On the other hand, one finds by direct (machine) computation that $\dim_{\gf}(\tansp_{[I]})$ $=$ $59$.  From this it follows that $[I]$ is a smooth point on an elementary component of $\Hilb^{13}_{\mathbb{A}^5_{\gf}}$ of dimension $59$; that is, $I$ is generic. Note that the dimension of the principal component in this case is $5 \cdot 13$ $=$ $65$.

\begin{rem} \label{rem:shapeRem}
Most of the examples presented in Section \ref{sec:examples} will be of the form just described: that is, there will be $n \geq 3$ variables $x_1, \dots, x_n$, with $x_1, \dots, x_{n-\kappa}$ the front variables and $x_{n-\kappa+1}, \dots, x_n$ the back variables ($1 < \kappa < n$).  The leading monomials will have degree $r$ $\geq$ $2$ and the trailing monomials will have degree $s$ $>$ $r$. The ideal $I$ will be generated by sufficiently general polynomials of the form
\[
  g_j\ =\ 
  \left( 
    \begin{array}{c}
     j\text{-th leading monomial}\ -\\
     \gf\text{-linear combination of trailing monomials}
    \end{array} \right).
\]
We say that ideals formed in this way have \textbf{shape} $(n,\kappa,r,s)$; they are special cases of a slightly more general type of ideal, introduced in Section \ref{sec:SLI}, that we call ``\SLI,'' and that is the focus of our exposition. Given an order ideal $\mathcal{O}$, one obtains a \SLI\ ideal by constructing its $\mathcal{O}$-border basis, making use of sets of leading and trailing monomials as in the preceding examples to obtain (some of) the generators. In particular, a \SLI\ ideal is in the $\mathcal{O}$-border basis scheme and is supported at one point (the origin) of $\mathbb{A}^n_{\gf}$.  The case of main interest (\SLI\ ideals of shape $(n,\kappa,r,s)$) is discussed in detail in Section \ref{sec:idealsOfShapeNkRs}.
\end{rem}

\subsection{Summary of examples} \label{ssec:sumOfEgs}

In Section \ref{sec:examples} we present several examples of generic ideals/elementary components, which we summarize briefly in the following list.  In each case, we give the Hilbert function, the shape (except for the last case), the dimension of the elementary component, which is equal to $\dim_{\gf}(\tansp_{[I]})$, and the dimension of the principal component.
 
\begin{description}
  \item[Hilbert func.\ (1,5,3,4,0), Shape $\mathbf{(5,2,2,3)}$] As discussed in Section \ref{ssec:AGen}, in this case $[I]$ is a smooth point on an elementary component of dimension $59$, and the dimension of the principal component is $5 \cdot 13 = 65$. 
  \item[Hilbert func.\ (1,5,3,4,5,6,0), Shape $\mathbf{(5,2,2,5)}$] In this case $[I]$ is a smooth point on an elementary component of dimension 104, and the dimension of the principal component is $5 \cdot 24$ $=$ $120$.
  \item[Hilbert func.\ (1,6,6,10,0), Shape $\mathbf{(6,3,2,3)}$] In this case $[I]$ is a smooth point on an elementary component of dimension $165$, and the dimension of the principal component is $6 \cdot 23$ $=$ $138$. 
  \item[Hilbert func.\ (1,6,21,10,15,0), Shape $\mathbf{(6,3,3,4)}$] In this case $[I]$ is a smooth point on an elementary component of dimension $705$, and the dimension of the principal component is $ 6 \cdot 53$ $=$ $318$.
  \item[Hilbert function (1,6,10,10,5,0)] In this case we give three different generic ideals/elementary components having this Hilbert function, each based on the same order ideal $\mathcal{O}$ of cardinality 32. The ideals have a slightly more general form than those in the examples above. In the first case, $[I]$ is smooth on an elementary component of dimension 255, in the second, the elementary component has dimension 222, and in the third, the dimension is 211. The principal component has dimension $6 \cdot 32$ $=$ $192$.  
\end{description}

\subsection{Plausible Genericity} \label{ssec:introPlausGen}

Unfortunately, we have no general theorems of the form ``every sufficiently general ideal of a certain type or shape is generic,'' since we do not know how to verify that $\dim_{\gf}(\tansp_{[I]})$ attains its minimum value except by direct computation in each case.  However, in Section \ref{sec:PlausArgs} we will make an attempt in this direction, by giving an easily-computable criterion for detecting if sufficiently general ideals of shape $(n,\kappa,r,s)$ are \emph{plausibly} generic, in which case we say that $(n,\kappa,r,s)$ is a \textbf{plausible} shape.  For example, this criterion indicates that the following shapes are plausible:
\begin{equation} \label{eqn:PlausibleShapes}
    \begin{array}{cc}
         \text{Shape}& \text{Range}\vspace{.05in}\\
         (n,2,2,3) & 5 \leq n \leq 50000 \text{ (at least)},\vspace{.02in}\\
         (n,2,2,4) & 5 \leq n \leq 50000\text{ (at least)},\vspace{.02in}\\
         (n,3,3,4) & 5 \leq n \leq 17,\vspace{.02in}\\
         (n,3,3,5) & 5 \leq n \leq 25,\vspace{.02in}\\
         (n,4,2,6) & 16 \leq n \leq 50000\text{ (at least)},\vspace{.02in}\\
         (n,4,3,6) & 8 \leq n \leq 120,\vspace{.02in}\\
         (50,\kappa,4,6) & 7 \leq \kappa \leq 22\vspace{.02in}\\
         (50,\kappa,4,8) &6 \leq \kappa \leq 14.
    \end{array}
\end{equation} 

Analysis of the asymptotic behavior of the plausibility criterion in Section \ref{ssec:PlausAsympBehavior} leads us to offer the following

\begin{conj} \label{conj:FirstOfTwo}
    Given $r=2$, $s>2$, and $\kappa \geq 2$, the shape $(n,\kappa,2,s)$  is plausible for all $n>>0$.  
\end{conj} 

\begin{rem} \label{rem:NoTrendRem}
    The reader may be wondering whether the trend suggested by the first two examples in Section \ref{ssec:sumOfEgs} extends.  The answer is almost certainly ``no.'' As $s$ increases, the shape 
$(5,2,2,s)$ must eventually become implausible, as shown in Section \ref{ssec:PlausAsympBehavior}.  Indeed, for the Hilbert function $(1,5,3,4,5,6,7,0)$, the tangent space dimension at $[I]$ is 139, which is less than 155, the dimension of the principal component, but is larger than the ``expected'' value of 131 (given by Equation (\ref{eqn:tanSpDimInGenShapeCase})) if $I$ were generic.  Hence, $[I]$ is not a point of the principal component, but neither is $I$ likely to be generic. 
\end{rem}

\subsection{Overview of paper} \label{ssec:overview}
Following the introduction, we review the terminology and theory of border basis schemes in Section \ref{sec:borderBases}, and lay some foundations for later sections.  In Section \ref{sec:SLI}, we present the definition and first properties of the ideals that are our main objects of study, which we call \textbf{\SLI\ ideals}.  (Their construction generalizes that described in Sections \ref{ssec:I-Eeg} and \ref{ssec:AGen}.)

Section \ref{sec:tanSpace} recalls the basic facts regarding the tangent space at $[I]$ for ideals $I$ of finite colength, and outlines how the dimension of the tangent space can be computed.  

Given a \SLI\ ideal $I$, in Section \ref{sec:idealFam} we  construct a map 
\[
    \fn{\EuScript{F}}{U}{\Hilb^{\mu}_{\mathbb{A}^n_{\gf}}},
\]
where $U$ is an affine space and $[I]$ $\in$ $\EuScript{F}(U)$.  Every point $[J]$ $\in$ $\EuScript{F}(U)$ corresponds to a closed subscheme $\Spec(\gf[\mathbf{x}]/J)$ supported at a single point of $\mathbb{A}^n_{\gf}$; in particular, $\EuScript{F}(U)$ contains all the points corresponding to ideals $I'$ obtained from $I$ through some combination of tweaking of coefficients and pulling back via automorphisms, as described in Section \ref{ssec:AGen}.  For $p$ $\in$ $U$ and $\EuScript{F}(p)$ $=$ $[I_p]$, we show (Proposition \ref{prop:F(U)DimLowerBd}) that the cardinality $L$ of a linearly independent set of vectors in the image of the derivative map  $\fn{\EuScript{F}'_p}{\tansp_p}{\tansp_{[I_p]}}$ is a lower bound  for the dimension of $\EuScript{F}(U)$. This leads to a simple method for finding elementary components (Proposition \ref{prop:genericCond}): If one has a lower bound $L$ for $\dim(\EuScript{F}(U))$ such that $\dim_{\gf}(\tansp_{[I_p]})$ $=$ $L$, then $[I_p]$ will be a smooth point on $\overline{\EuScript{F}(U)}$, which is accordingly an elementary component of dimension $L$.

We study the derivative map $\fn{\EuScript{F}'_p}{\tansp_p}{\tansp_{[I_p]}}$ in Section \ref{sec:derivMap}.  In particular, we compute the images of the standard unit vectors at $p \in U$ (an affine space) to enable us to find lower bounds $L$ on the dimension of $\EuScript{F}(U)$.   

In Section \ref{sec:spCaseRestr} we define and study the \textbf{lex-segment complement} order ideals (and associated \SLI\ ideals) that are used in all of the examples outlined in Section \ref{ssec:sumOfEgs}.  We work out in detail the concepts and results of Sections \ref{sec:idealFam} and \ref{sec:derivMap} in this special case, to prepare for the presentation  of the examples in Section \ref{sec:examples}.

Following the presentation of the examples, we turn to the final goal of the paper, which is to develop a criterion for detecting plausible shapes $(n,\kappa,r,s)$.  The criterion is stated and justified in (the final) Section \ref{sec:PlausArgs}, preceded by an extensive preparatory study of \SLI\ ideals of shape $(n,\kappa,r,s)$ in Section \ref{sec:idealsOfShapeNkRs}.

\medskip

\noindent \textit{Acknowledgment}: The author is happy to thank Prof.\ Tony Iarrobino for several helpful private communications.  He also thanks the referee for a myriad of useful comments.


\section{Border basis schemes} \label{sec:borderBases}

In this section, we briefly recall some of the terminology and theory of \textbf{border basis schemes} as given in  \cite[Secs.\ 2, 3]{KreutzerAndRobbiano1:DefsOfBorderBases}.

\subsection{Basic definitions}  \label{ssec:BasDefs} One begins with an \textbf{order ideal}, which is a finite set $\mathcal{O}$ $=$ $\{t_1,\dots t_{\mu} \}$ of monomials in the variables $x_1,\, \dots,\,  x_n$ such that whenever a monomial $m$ divides a member of $\mathcal{O}$, it follows that $m \in \mathcal{O}$.  We will refer to the monomials $t_i$ as \textbf{basis monomials}. The \textbf{border} of $\mathcal{O}$ is the set of monomials 
\[
  \partial \mathcal{O}\ =\  \left(x_1 \mathcal{O} \cup \dots \cup x_n \mathcal{O}\right) \setminus \mathcal{O} \ =\ \{b_1,\, \dots,\, b_{\nu}\};
\]
we will refer to the $b_j$ as \textbf{boundary monomials}.  A set of polynomials $\mathcal{B}$ $=$ $\{ g_1,\dots, g_{\nu} \}$ of the form $g_j = b_j - \sum_{i = 1}^{\mu}c_{ij}t_i$ with $c_{ij} \in \gf$ is called an \textbf{$\mathcal{O}$-border prebasis} of the ideal
\[
   I = (\mathcal{B}) \subseteq \gf[x_1,\, \dots,\, x_n] =  \gf[\mathbf{x}] = R.
\]
It is clear that every boundary monomial is congruent to a linear combination of basis monomials modulo $I$, and an induction argument shows that the same is true for every monomial; in other words, the quotient $R/I$ is spanned as a $\gf$-vector space by the $t_i$.  We say that $\mathcal{B}$ is an \textbf{$\mathcal{O}$-border basis} of $I$ if $\mathcal{O}$ is a $\gf$-basis of the quotient; in this case, every monomial is congruent modulo $I$ to a unique $\gf$-linear combination of basis monomials. 

The $\mathcal{O}$\textbf{-border basis scheme} $\mathbb{B}_{\mathcal{O}}$ is an affine scheme whose $\gf$-points correspond to the ideals $I$ having an $\mathcal{O}$-border basis; as such, it is an open affine subscheme of $\Hilb^{\mu}_{A^n_{\gf}}$, from which it inherits the universal property (\ref{txt:BbasUnivProp}): Let $\mathcal{Z}_{\mathcal{O}}$ denote the restriction of the universal closed subscheme $\mathcal{Z}_{\mu}$ $\subseteq$ $\Hilb^{\mu}_{\mathbb{A}^n_{\gf}} \times \mathbb{A}^n_{\gf}$ to $\mathbb{B}_{\mathcal{O}}$.  Then (see, \eg, \cite[Th.\ 37, p.\ 306]{Huib:UConstr}):

\numbtext{txt:BbasUnivProp}{.}{.}{A map $\fn{q}{Q}{\mathbb{B}_{\mathcal{O}}}$ corresponds uniquely to a closed subscheme $\mathcal{Z}_q$ $\subseteq$ $Q \times \mathbb{A}^n_{\gf}$ such that the direct image of $O_{Z_q}$ on 
$Q$ is free with basis $\mathcal{O}$; the correspondence is given by $q$ $\leftrightarrow$ $q^*(\mathcal{Z}_{\mathcal{O}})$.}

Note that the schemes $\mathbb{B}_{\mathcal{O}}$, as $\mathcal{O}$ ranges over all order ideals of cardinality $\mu$, form an open affine covering of $\Hilb^{\mu}_{\mathbb{A}^n_{\gf}}$; moreover, one can construct  $\Hilb^{\mu}_{\mathbb{A}^n_{\gf}}$ by first constructing the schemes $\mathbb{B}_{\mathcal{O}}$ and then gluing along the natural overlaps (see, \eg, \cite{LaksovAndSkjelnes:HilbSchConstr}, \cite{Huib:UConstr},
 \cite{KreuzerAndRobbiano2:TheGeometryOfBorderBases}).


\subsection{Neighbor syzygies of a border basis} \label{ssec:NbrSyz} Let $I$ be an ideal having $\mathcal{O}$-border basis $\mathcal{B}$ $=$ $\{ g_1, \dots, g_{\nu} \}$ $\subseteq$ $\gf[\mathbf{x}]$, as in Section \ref{ssec:BasDefs}. The neighbor syzygies provide a convenient set of generators for the first syzygy module of the $g_j$; we briefly recall their construction, following \cite[Sec.\ $4$, p.\ $13$]{KreutzerAndRobbiano1:DefsOfBorderBases}.

We say that two boundary monomials $b_j$ and $b_{j'}$ are \textbf{next-door neighbors} if $x_k b_j$ $=$ $b_{j'}$ for some $k$ $\in$ $\{ 1,\dots,n \}$, and \textbf{across-the-street} neighbors if $x_k b_j$ $=$ $x_l b_{j'}$ for some $k,\l$ $\in$ $\{ 1,\dots,n \}$; in either case, we say that $b_j$ and $b_{j'}$ are \textbf{neighbors}.  Given neighbors $b_j$ and $b_{j'}$, we form the $S$-polynomial   
\[
     S(g_j,g_{j'}) = \left\{ \begin{array}{ll}
                                  x_k \cdot g_j - g_{j'}, & \text{ if } b_j, b_{j'} \text{ are next-door nbrs,}\vspace{.05in}\\
                                  x_k \cdot g_j - x_l \cdot g_{j'}, & \text{ if } b_j, b_{j'} \text{ are across-the-street nbrs.} 
                             \end{array} \right. 
\]
In either case, $S(g_j,g_{j'})$ is a $\gf$-linear combination of basis and boundary monomials.  For each term of the form $c_{j''} b_{j''}$ that appears in $S(g_j,g_{j'})$, we subtract $c_{j''} g_{j''}$ (put another way, we reduce $S(g_j,g_{j'})$ modulo the ideal generators $\mathcal{B}$); the result is a $\gf$-linear combination of basis monomials $\overline{S}(g_j,g_{j'})$ that is a $\gf[\mathbf{x}]$-linear combination of the $g_j$; whence, $\overline{S}(g_j,g_{j'})$ $\equiv 0 {\mod I}$.  Since $\mathcal{B}$ is an $\mathcal{O}$-border basis of $I$, it follows at once that $\overline{S}(g_j,g_{j'})$ is the $0$-polynomial, so we have constructed a syzygy of the polynomials $g_j$, the \textbf{neighbor syzygy} associated to the neighbors $b_j$, $b_{j'}$ $\in$ $\partial \mathcal{O}$.  Writing this syzygy in the form $\sum_{\hat{j}=1}^{\nu}f_{\hat{j}}\, g_{\hat{j}}$ $=$ $0$, the tuple of coefficients $(f_{\hat{j}})$ has at most two components of degree 1 ($f_j$ and possibly $f_{j'}$), and the remaining components are all constants.  

If $\mathcal{B}$ $=$ $(g_1, \dots, g_{\nu})$ is just an $\mathcal{O}$-border prebasis, one can still compute the $S$-polynomials and their reductions $\overline{S}(g_j,g_{j'})$.  These again are $\gf$-linear combinations of basis monomials, but they no longer necessarily vanish.  We have the following key results:

\begin{prop} \label{prop:NbrSyzFacts} {  }
\noindent

\begin{itemize}
  \item[i.]    
 If the reduced $S$-polynomials $\overline{S}(g_j,g_{j'})$ for an $\mathcal{O}$-border prebasis $\mathcal{B}$ are all equal to $0$, then $\mathcal{B}$ is an $\mathcal{O}$-border basis; that is, the quotient $\gf[\mathbf{x}]/(\mathcal{B})$ is $\gf$-free with basis $\mathcal{O}$.  (The converse was proved at the start of this section.)
  \item[ii.]  If $\mathcal{B}$ is an $\mathcal{O}$-border basis, then the neighbor syzygies generate the first syzygy module of the $\{g_j\}$ as a $\gf[\mathbf{x}]$-module.
\end{itemize} 
\end{prop}

\proof The first statement is proved in \cite[Prop.\ 6.4.34, p.\ 438]{KreutzerAndRobbianoVolTwo}, and both statements are proved in  \cite[Th.\ $22$, p.\ 292]{Huib:UConstr}. A beautiful algorithmic proof of the second statement is given in \cite{KreuzerAndKriegl}.
\qed

\medskip

\begin{rem} \label{rem:GrdRingsOK}
     The theory of border bases can be developed in essentially the same way as summarized above when the ground field $\gf$ is replaced by an arbitrary commutative and unitary ring $A$, such as a $\gf$-algebra (see, \eg, \cite{Huib:UConstr}).  In particular, the analogue of Proposition  \ref{prop:NbrSyzFacts} holds in this more general context.
\end{rem}


\subsection{Linear syzygies} \label{ssec:AlgorithmForNbrSyz}

We say that a syzygy $(f_j)$ of the ideal generators $g_j$ of an $\mathcal{O}$-border basis (that is, $\sum_{h=1}^{\nu}(f_j\cdot g_j) = 0$) is a \textbf{linear syzygy} provided that the coefficients $f_j$ $\in$ $\gf[\mathbf{x}]$ have degree at most 1.  For example, the neighbor syzygies are all linear syzygies.  Since the neighbor syzygies generate the $\gf[\mathbf{x}]$-module of first syzygies of the border basis $\mathcal{B}$, by Proposition \ref{prop:NbrSyzFacts}, we see that a $\gf$-basis of the linear syzygies is also a set of $\gf[\mathbf{x}]$-generators of the full syzygy module.

We briefly describe the algorithm we use for computing a $\gf$-basis of the linear syzygies; as a consequence, we will obtain the cardinality of this basis.  (This algorithm is implemented in the \emph{Mathematica} function \textbf{makeLinearSyzygies}, included in the notebook of utility functions mentioned at the start of Section \ref{sec:examples}.)  

We first compute the set of boundary monomials $\partial \mathcal{O}$ and the set of \textbf{target monomials} 
\begin{equation}\label{eqn:tarMonDef}
    \tarmons = \partial \mathcal{O} \cup \{ x_{\alpha} \cdot b_j \mid 1 \leq \alpha \leq n,\ b_j \in \partial \mathcal{O}\}.
\end{equation}
Next we define  a $\gf$-linear ``projection'' map $\fn{\pi_T}{\gf[\mathbf{x}]}{\Span_{\gf}(\tarmons)}$ by extending linearly the map on monomials 
\[
    m \mapsto \left\{ \begin{array}{l}
                           m,\ \text{if } m \in \tarmons,\vspace{.05in}\\
                           0,\ \text{otherwise}.
                     \end{array} \right.
\]  
Let $V$ be a $\gf$-vector space with basis
\[
    E\ =\ \{ e_{\alpha,j} \mid 0 \leq \alpha \leq n,\ 1 \leq j \leq \nu\},\ \text{ of cardinality } |E| = (n+1)\cdot \nu,
\]
and define a linear map
\[
    \fn{\sigma}{V \approx \gf^{(n+1)\cdot \nu}}{\Span_{\gf}(\tarmons)},\ \  e_{\alpha,j} \mapsto 
    \left\{ \begin{array}{l}
        \pi_{\tarmons}(1\cdot g_{j}) = b_j,\ \text{if } \alpha = 0,\vspace{.05in}\\
        \pi_{\tarmons}(x_{\alpha} \cdot g_{j}),\ \text{if } n \geq \alpha > 0
    \end{array} \right. .
\]

\begin{lem} \label{lem:sigmaSurjLem}
    The map $\sigma$ is surjective.
\end{lem}

\proof
    Since 
\[
    e_{1,j} \mapsto b_j \in \partial \mathcal{O} \subseteq \tarmons \text{ for } 1 \leq j \leq \nu,
\]
it is clear that $\Span_{\gf}(\partial \mathcal{O})$ is in the image of $\sigma$.  We now observe that for $0< \alpha < n$, 
\[
    e_{\alpha,j} \mapsto \pi_{\tarmons}(x_{\alpha} \cdot g_j) = x_{\alpha} \cdot b_j + \left( \text{element of } \Span_{\gf}(\partial \mathcal{O})\right);
\]
whence, every monomial $x_{\alpha} \cdot b_j$ $\in$ $\tarmons$ is in the image of $\sigma$, and the lemma follows at once. 
\qed

\medskip

Let $(d_{\alpha,j})$ $\in$ $\operatorname{Ker}(\sigma)$. Setting $f_j$ $=$ $d_{0,j} + \sum_{\alpha=1}^n d_{\alpha,j} x_{\alpha}$, we observe that the tuple $(f_j)$ is a linear syzygy of $\{g_j\}$, because $\sum_{j=1}^{\nu}{f_j}\cdot{g_j}$ is a $\gf$-linear combination of basis monomials (the monomials in $\tarmons$ having cancelled out), and hence $0$. Moreover, every linear syzygy of the $g_j$ arises in this way.  So a $\gf$-basis of the linear syzygies can be computed simply by computing a basis $\{ (d_{\alpha,j}) \}$ of the kernel of $\sigma$ and assembling the corresponding tuples $(f_j)$. From this it follows that the dimension of the $\gf$-vector space of linear syzygies of the border  basis $\mathcal{B}$ of $I$ is given by
\begin{equation} \label{eqn:linSyzCount}
   \psi\ =\ \dim_{\gf}(\ker(\sigma))\ =\ (n+1) \cdot \nu - |\tarmons|.
\end{equation}


\subsection{Generators of the ideal of $\mathbb{B}_{\mathcal{O}}$} \label{ssec:BorBasGens}  The border basis scheme 
 $\mathbb{B}_{\mathcal{O}}$ is a closed subscheme of $\mathbb{A}^{\mu \nu}_{\gf}$ $=$ $\Spec(\gf[\mathcal{C}])$, where 
\begin{equation} \label{eqn:setCDefn}
    \mathcal{C} = \{C_{ij},\ 1 \leq i \leq \mu,\ 1 \leq j \leq \nu\}
\end{equation}
 is a set of indeterminates; the point corresponding to the ideal $I$ having border basis $\{g_j = b_j - \sum_{i = 1}^{\mu}c_{ij}t_i\}$ is $(c_{ij})$ $\in$ $\mathbb{A}^{\mu \nu}_{\gf}$.  The generators of the ideal $\mathcal{I}_{\mathcal{O}}$ such that  $\mathbb{B}_{\mathcal{O}}$ $=$ $\Spec(\gf[C_{ij}])/\mathcal{I}_{\mathcal{O}}$ can be obtained as follows: Form the ``generic $\mathcal{O}$-border prebasis'' 
\begin{equation} \label{eqn:genBorPreBasis}
    \mathcal{B}^{\star} = \{ G_1, \dots, G_{\nu} \} \subseteq \gf[\mathcal{C}][\mathbf{x}],\ \ G_j = b_j - \sum_{i=1}^{\mu}C_{ij}t_i,
\end{equation}
and compute the $\gf[\mathcal{C}]$-linear combinations of basis monomials 
\[
        \overline{S}(G_j,G_{j'}) = \sum_{i=1}^{\mu} \varphi^{j,j'}_i t_i.
\] 
Then the ideal $\mathcal{I}_{\mathcal{O}}$ is generated by the coefficients $\varphi^{j,j'}_i$ $\in$ $\gf[\mathcal{C}]$ (see, \eg, \cite[Th.\ 37, p.\ 306]{Huib:UConstr}.  The point is that over the ring $A_{\mathcal{O}}$ $=$ $\gf[\mathcal{C}]/\mathcal{I}_{\mathcal{O}}$, the polynomials $\overline{S}(G_j,G_{j'})$ all vanish, so Proposition \ref{prop:NbrSyzFacts} (in the light of Remark \ref{rem:GrdRingsOK}) yields that the quotient $A_{\mathcal{O}}[\mathbf{x}]/(G_j)$ is $A_{\mathcal{O}}$-free with basis $\mathcal{O}$.  One then shows that $\Spec(A_{\mathcal{O}})$ and the family of subschemes $\Spec(A_{\mathcal{O}}[\mathbf{x}]/(G_j))$ together satisfy the universal property (\ref{txt:BbasUnivProp}) of the border basis scheme $\mathbb{B}_{\mathcal{O}}$.

\begin{rem} \label{rem:matrixMethod}
    The generators of the ideal $\mathcal{I}_{\mathcal{O}}$ can also be constructed as the entries of the commutators of the ``generic multiplication matrices'' --- see, \eg, \cite[Sec.\ 3]{KreutzerAndRobbiano1:DefsOfBorderBases}.
\end{rem}



\section{\SLICap\ ideals} \label{sec:SLI}

In this section we will describe our main objects of study. 


\subsection{Definition of \SLI\ ideals} \label{ssec:DistIdealsDefn}

Let $\mathcal{O}_{\text{max}}$ $\subseteq$ $\mathcal{O}$ denote the subset of  \textbf{maximal} basis monomials, which are those basis monomials $t_i$ such that  $x_k t_{i}$ $\in$ $\partial \mathcal{O}$ for all $1 \leq k \leq n$.  Similarly, let $\partial \mathcal{O}_{\text{min}}$ $\subseteq$ $\partial \mathcal{O}$ denote the subset of \textbf{minimal} boundary monomials, which are those boundary monomials $b_j$ such that $b_j/x_k$ $\in$ $\mathcal{O}$ for every $x_k$ that appears in $b_j$.
Choose non-empty subsets
\[
  \leadmon = \{b_{j_1},\, \dots,\, b_{j_{\lambda}}\} \subseteq \partial \mathcal{O}_{\text{min}}
\]
and
\[
     \trailmon = \{t_{i_1},\, \dots,\, t_{i_{\tau}} \} \subseteq \mathcal{O}_{\text{max}} 
\]
such that the set $\leadmon$ is disjoint from the subset of boundary monomials 
\begin{equation} \label{eqn:bdryT}
    \partial \trailmon = \{x_k t_{i_{\ell}} \mid 1 \leq k \leq n,\ t_{i_{\ell}} \in \trailmon \} \subseteq \partial \mathcal{O}.
\end{equation}
We will call $\leadmon$ (resp.\ $\trailmon$) the \textbf{leading} (resp.\ \textbf{trailing}) \textbf{monomials}.

We choose a set $G$ $=$ $\{ g_{j_{\iota}} \}$ of polynomials of the form
\begin{equation} \label{eqn:fDef}
	g_{j_{\iota}} = b_{j_{\iota}} - N_{j_{\iota}},\ \ 1 \leq \iota \leq \lambda,\ b_{j_{\iota}} \in \leadmon,\ N_{j_{\iota}} \in {\rm Span}_{\gf}(\trailmon),
\end{equation}
and extend $G$ to an $\mathcal{O}$-border prebasis 
\[
    \mathcal{B} = \{ g_1, \dots, g_{\nu} \} = G \cup \left( \partial \mathcal{O} \setminus \leadmon \right).
\]
In Proposition \ref{prop:specialCCor} we prove that $\mathcal{B}$ is an $\mathcal{O}$-border basis of an ideal $I$ $=$ $(\mathcal{B})$.  We say that any such ideal $I$ is a \textbf{\SLI\ ideal}.

\subsection{Example} \label{ssec:RunningExample1}
For the order ideal 
\[
    \mathcal{O} =  \left\{1,\, x_1,\, x_2,\, x_3,\, x_2 x_3,\, x_3^2\right\} \subseteq \gf[x_1,x_2,x_3], \text{ with } 
\mu = |\mathcal{O}| = 6,
\]
one has that
\[
    \begin{array}{rcl}
        \partial \mathcal{O} & = & \left\{
        \begin{array}{c} x_1^2,\ x_1\, x_2,\ x_1\, x_3,\ x_2^2,\vspace{.05in}\\
            x_1\, x_2\, x_3,\ x_1\, x_3^2,\ x_2^2\, x_3,\ x_2\, x_3^2,\ x_3^3
        \end{array} \right\},\vspace{.05in}\\
        \mathcal{O}_{\text{max}} & = & \left\{x_1,\, x_2 x_3,\, x_3^2\right\},\vspace{.05in}\\
        \partial \mathcal{O}_{\text{min}} & = &\left\{x_1^2,\, x_1 x_2,\, x_1 x_3,\, x_2^2,\, x_2 x_3^2,\, x_3^3\right\}.
    \end{array}
\]
There are various possible choices for the sets $\leadmon$ and $\trailmon$ of leading and trailing monomials that satisfy $\leadmon \cap \partial \trailmon$ $=$ $\emptyset$; here is one:
\[
         \leadmon = \{ x_1^2,\, x_1 x_2,\, x_1 x_3,\, x_2^2 \},\  \trailmon = \{ x_2 x_3,\, x_3^2 \},
\]
so $\lambda =  |\leadmon| =  4$ and $\tau = |\trailmon| = 2$.  We therefore have an ($8 = \lambda \mu$)-dimensional family of \SLI\ ideals in the border basis scheme $\mathbb{B}_{\mathcal{O}}$ $\subseteq$ $\Hilb^6_{\mathbb{A}^3_{\gf}}$ with border bases of the form
\begin{equation} \label{eqn:RunningExampleBorBases}
  \mathcal{B}\ =\  \left\{ \begin{array}{l}
             g_1 = x_1^2 - C_{5,1}\,x_2\,x_3 - C_{6,1}\,x_3^2,\vspace{.05in}\\
             g_2 = \ x_1\, x_2 - C_{5,2}\,x_2\,x_3 - C_{6,2}\,x_3^2,\vspace{.05in}\\
             g_3 = x_1\, x_3 - C_{5,3}\,x_2\,x_3 - C_{6,3}\,x_3^2,\vspace{.05in}\\
             g_4 = x_2^2 - C_{5,4}\,x_2\,x_3 - C_{6,4}\,x_3^2,\vspace{.05in}\\
             g_5 = x_1\, x_2\, x_3,\ \ g_6 =  x_1\, x_3^2,\ \ g_7 =  x_2^2\, x_3,\vspace{.05in}\\ 
             g_8 =  x_2\, x_3^2,\ \ g_9 = x_3^3             
    \end{array} \right\}
\end{equation}
(the indeterminate coefficients $C_{ij}$ $\in$ $\mathcal{C}$ (\ref{eqn:setCDefn}) would of course be replaced by elements of $\gf$ in any specific example).  To verify that $\mathcal{B}$ is an $\mathcal{O}$-border basis, it suffices, by Proposition \ref{prop:NbrSyzFacts}, to show that all the reduced $S$-polynomials $\overline{S}(g_j,g_{j'})$ are equal to $0$.  The general argument is given in the proof of Proposition \ref{prop:specialCCor}; we can easily check this by hand for the pre-basis $\mathcal{B}$ in (\ref{eqn:RunningExampleBorBases}); here, for example, is one of the required verifications:  
\[
    \begin{array}{rcl}
       S(g_1,g_2) & = & x_2 \cdot g_1 - x_1\cdot g_2\vspace{.05in}\\
           {}     & = & x_2 (x_1^2 - C_{5,1}\,x_2\,x_3 - C_{6,1}\,x_3^2)\ - \vspace{.05in}\\
           {}     &{} & x_1(x_1\, x_2 - C_{5,2}\,x_2\,x_3 - C_{6,2}\,x_3^2)\vspace{.05in}\\
           {}     & = & - C_{5,1}\,x_2^2\,x_3 - C_{6,1}\,x_2\,x_3^2 + C_{5,2}\,x_1\,x_2\,x_3 + C_{6,2}\,x_1\,x_3^2\vspace{.05in}\\
           {}     & \longrightarrow_{\mathcal{B}} & 0.
    \end{array}    
\]

\subsection{The locus of \SLI\ ideals in $\mathbb{B}_{\mathcal{O}}$} \label{ssec:specialLoci}

We define a subset $\specl$ of $\mathcal{C}$ (\ref{eqn:setCDefn}) as follows:
\begin{equation}\label{eqn:SpecialC}
    \specl\ =\ \{ C_{ij} \mid  b_j \in \leadmon,\ t_i \in \trailmon\}\ \subseteq\ \mathcal{C}.
\end{equation}
We will say that the members of $\specl$ and the associated pairs of indices $(i,j)$ are \textbf{\SLI}.

Consider the surjection of polynomial rings

\[
   \fn{\gamma}{\gf[\mathcal{C}]}{\gf[\specl]}, \ \ C_{ij} \mapsto \left\{ 
    \begin{array}{l}
        C_{ij}, \text{ if } C_{ij} \in \specl,\vspace{.05in}\\
        0, \text{ otherwise},
    \end{array} \right.,
\]
and let $\fn{\hat{\gamma}}{\gf[\mathcal{C}][\mathbf{x}]}{\gf[\specl}][\mathbf{x}]$ denote the map obtained by applying $\gamma$ to each coefficient of the input polynomial $f \in \gf[\mathcal{C}][\mathbf{x}]$).

\begin{prop} \label{prop:specialCCor}
The map $\gamma$ factors through the coordinate ring $\gf[\mathcal{C}]/\mathcal{I}_{\mathcal{O}}$ of $\mathbb{B}_{\mathcal{O}}$.  Consequently, $\mathbb{B}_{\mathcal{O}}$ contains an $|\specl| = \lambda \tau$-dimensional closed subscheme $X_{\specl}$ $=$ $\Spec(\gf[\specl])$ isomorphic to affine space, and whose $\gf$-points are obtained by assigning arbitrary scalars to the indeterminates in $\specl$ and $0$ to the other indeterminates. Moreover, every point $[I']$ $\in$ $X_{\specl}$ corresponds to a closed subscheme $\Spec(\gf[\mathbf{x}]/I')$ that is supported at the origin of $\mathbb{A}^n_{\gf}$.
\end{prop}

\proof
      The image under $\hat{\gamma}$ of the $\mathcal{O}$-border prebasis $\mathcal{B}^{\star}$ 
(\ref{eqn:genBorPreBasis}) has the form
\[
    \{\hat{\gamma}(G_1),\, \dots,\, \hat{\gamma}{G_{\nu}} \},\ \hat{\gamma}(G_j) = \left\{ 
                                             \begin{array}{l}
                                                 b_j, \text{ if } b_j \notin \leadmon,\vspace{.05in}\\
                                                 b_j - \sum_{C_{ij} \in \specl} C_{ij}t_i, \text{ if } 
                                                 b_j  \in \leadmon
                                             \end{array} \right. .
\]
We claim that the polynomials $\overline{S}(\hat{\gamma}(G_j),\hat{\gamma}(G_{j'}))$ of Section \ref{ssec:NbrSyz} (where the reductions are with respect to $\hat{\gamma}(\mathcal{B}^{\star})$) all vanish.  There are three cases to check: First suppose that $b_j$ and $b_{j'}$ are neighbors such that neither is a leading monomial.  Then 
\[
    S(\hat{\gamma}(G_j),\hat{\gamma}(G_{j'})) = S(b_j,b_{j'}) = 0\ \Rightarrow\ \overline{S}(\hat{\gamma}(G_j),\hat{\gamma}(G_{j'})) = 0.
\] 
The second case to consider is that of two neighbors $b_j$ and $b_{j'}$ such that $b_j$ $\in$ $\leadmon$ and $b_{j'}$ $\notin$ $\leadmon$.  Note that we cannot have that $b_j$ $=$ $x_j\, b_{j'}$ because $b_j$ is a minimal boundary monomial.  It follows that
\[
  \begin{array}{rcl}
    S(\hat{\gamma}(G_j),\hat{\gamma}(G_{j'})) &  = &  x_k \left( b_j- \sum_{C_{ij} \in \specl} C_{ij}t_i \right) - (x_l\text{ or } 1)\, b_{j'}\vspace{.05in}\\
     {} & {=} & - \left( \sum_{C_{ij} \in \specl} C_{ij}\, (x_k\, t_i) \right).
  \end{array}
\]
Since the only terms that survive in the last expression have \SLI\ coefficients $C_{ij}$, we know that $t_i$ $\in$ $\trailmon$ and therefore 
\[
    x_k\, t_i \in \partial \mathcal{O} \setminus \leadmon\ \Rightarrow x_k\,t_i \in \hat{\gamma}(\mathcal{B}^{\star})\ \Rightarrow\ \overline{S}(\hat{\gamma}(G_j),\hat{\gamma}(G_{j'})) = 0.
\] 
The third case is that of neighbors $b_j$, $b_{j'}$ $\in$ $\leadmon$, for which the argument is similar to that of the second case (and is illustrated in the Example of Section \ref{ssec:RunningExample1}).

  Since it is clear that 
\[
  \begin{array}{rcl}
0\ =\ \overline{S}(\hat{\gamma}(G_j)),\hat{\gamma}(G_{j'})) & = & \hat{\gamma}(\overline{S}(G_j,G_{j'}))\vspace{.05in}\\  
{} & = & \hat{\gamma} \left( \sum_{i=1}^{\mu} \varphi^{j,j'}_i t_i \right)\ \vspace{.05in}\\ 
{} & = &   \sum_{i=1}^{\mu} \gamma(\varphi^{j,j'}_i) t_i ,
  \end{array}
\]
we conclude that $\gamma(\varphi^{j,j'}_{i})$ $=$ $0$ for all neighbor pairs $b_j$, $b_{j'}$ and all $1 \leq i \leq \mu$.  This proves the first part of the proposition.

To prove the last statement, it suffices to show that for each variable $x_k$, there is an exponent $e_k$ such that $x_k^{e_k}$ $\in$ $I'$.  To this end, let $e'_k$ be the least $e$ such that $x_k^e$ $\notin$ $\mathcal{O}$, in which case $x_k^{e'_k}$ $=$ $b_{j_k}$ $\in$ $\partial \mathcal{O}$.  If $b_{j_k}$ $\notin$ $\leadmon$, then $I'$ contains the polynomial $g_{j_k}$ $=$ $b_{j_k}$ (indeed, $I'$ contains every boundary monomial $b_j$ $\notin$ $\leadmon$).  On the other hand, if $b_{j_k}$ $\in$ $\leadmon$, then $I'$ contains the polynomial $g_{j_k}$ $=$ $b_{j_k}- \sum_{i=1}^{\mu} c_{i j_k}t_i$, for which the coefficient $c_{i,j_k}$ $\neq$ $0$ $\Rightarrow$ $t_i$ $\in$ $\trailmon$.  Multiplying this polynomial by $x_k$, and recalling that  $t_i$ $\in$ $\trailmon$ $\Rightarrow$ $x_k \cdot t_i$ $\in$ $\partial \mathcal{O}$ $\setminus$ $\leadmon$,  we see that $I'$ contains a polynomial of the form
\[
    x_k \cdot b_{j_k} - (\gf\text{-linear combination of monomials in } \partial \mathcal{O}\setminus \leadmon),
\]
which implies that $x_k \cdot b_{j_k}$ $=$ $x_k^{e'_k + 1}$ $\in$ $I'$, thereby completing the proof. \qed

\begin{rem} \label{rem:specialCRem}
    A special case of this proposition appeared in \cite[Cor.\ 41, p.\ 313]{Huib:UConstr}.
\end{rem}

Evidently the $\gf$-points of $X_{\specl}$ are the points $[I]$ $\in$ $\Hilb^{\mu}_{\mathbb{A}^n_{\gf}}$ corresponding to the \SLI\ ideals $I$. The proposition then immediately yields the following

\begin{cor} \label{cor:DistIdealsSuppAtOnePt}
    Every \SLI\ ideal is supported at a single point (the origin) of $\mathbb{A}^n_{\gf}$.  
\qed
\end{cor}

 We will call the locus $X_{\specl}$ $\cong$ $\mathbb{A}^{\lambda \tau}_{\gf}$ the \textbf{\SLI\  locus} associated to $\mathcal{O}$, $\leadmon$, and $\trailmon$.  The monomial ideal $I_0$ $=$ $(\partial \mathcal{O})$ corresponds to  the origin of $X_{\specl}$, that is, the point $[I_0]$ $\in$ $X_{\specl}$ defined by setting all the \SLI\ indeterminates $C_{ij}$ to $0$.

%
%

\subsection{Efficient \SLI\ ideals} \label{ssec:EffSLIIdeals}

   We say that a \SLI\ ideal $I$ is \textbf{efficient} provided that $I$ $=$ $(\mathcal{B})$ is generated by the subset $G$ $=$ $\{ g_{j_{\iota}} \}$ $\subseteq$ $\mathcal{B}$ (equation (\ref{eqn:fDef})).  Since in any case $(G)$ $\subseteq$ $(\mathcal{B})$, we have that $(G)$ $=$ $(\mathcal{B})$ if and only if $(\mathcal{B})$ $\subseteq$ $(G)$, which is equivalent to $\partial \mathcal{O} \setminus \leadmon$ $\subseteq$ $(G)$. This property is easy to test computationally: Simply compute a Groebner basis for $(G)$ and reduce each non-leading boundary monomial $b_j$  modulo the Groebner basis; $I$ is efficient if and only if all the reductions are $0$.

\subsection{Example (continued)} \label{ssec:RunningExampleII}

Continuing with the Example of Section \ref{ssec:RunningExample1}, suppose that all of the coefficients $C_{ij}$ in (\ref{eqn:RunningExampleBorBases}) are set to $0$.  Then
\[
  \begin{array}{rcl}
    G           & = &\{ x_1^2,\ x_1\, x_2,\ x_1\, x_3,\ x_2^2 \}, \text{ and}\vspace{.05in}\\
    \mathcal{B} & = & G \cup \{ x_1\, x_2\, x_3,\ x_1\, x_3^2,\ x_2^2\, x_3,\ x_2\, x_3^2,\ x_3^3 \}.
  \end{array}
\]
Since, for example, $x_3^3$ $\notin$ $(G)$, we have that $(G)$ $\neq$ $(\mathcal{B})$, so $(\mathcal{B})$ is not efficient in this case.  On the other hand, if we take
\[
  \begin{array}{rcl}
    G &  = & \left\{x_1^2 +  x_2\, x_3 + x_3^2,\ x_1\, x_2 + x_2\, x_3,\ x_1\, x_3 +x_2\, x_3 + x_3^2,\ x_2^2\right\}, \text{ and}\vspace{.05in}\\
    \mathcal{B} & = & G \cup  \{ x_1\, x_2\, x_3,\ x_1\, x_3^2,\ x_2^2\, x_3,\ x_2\, x_3^2,\ x_3^3 \},
  \end{array}
\]
then we have that $(\mathcal{B})$ is efficient. To show this, we first compute a (lex) Groebner basis of $(G)$; the result is
\[
    \left\{x_1^2 +  x_2\, x_3 + x_3^2,\ x_1\, x_2 + x_2\, x_3,\ x_1\, x_3 +x_2\, x_3 + x_3^2,\ x_2^2,\ 
        x_2\, x_3^2,\ x_3^3 \right\},  
\]
which shows that $x_2\, x_3^2,\ x_3^3 \in (G)$.  It follows that $x_3 \cdot \trailmon$ $\subseteq$ $(G)$; since in addition $x_3 \cdot G$ $\subseteq$ $(G)$, we obtain that $x_3 \cdot \leadmon$ $\subseteq$ $(G)$; whence, $\{x_1\,x_2\,x_3,\, x_1\,x_3^2,\, x_2^2\, x_3 \}$ $\subseteq$ $(G)$, so $(\mathcal{B})$ $\subseteq$ $G$, and we are done.   

\begin{rem} \label{rem:RunningExampleNotGeneric}
    We noted in the Introduction that a non-trivial elementary component can only exist for $\mu \geq 8$; therefore, none of the \SLI\ ideals associated to $\mathcal{O}$, $\leadmon$, and $\trailmon$ as in Section \ref{ssec:RunningExample1} can be generic.  Indeed, for the efficient \SLI\ ideal $(\mathcal{B})$ just discussed, the tangent space dimension at $[(\mathcal{B})]$ is $18$, so that $[(\mathcal{B})]$ is a smooth point on the irreducible variety $\Hilb^{6}_{\mathbb{A}^3_{\gf}}$.
\end{rem}


\subsection{A sufficient condition for efficiency} \label{ssec:effTesting}

By analogy with $\partial \trailmon$ (\ref{eqn:bdryT}), we let
\[
    \partial \leadmon = \{x_k b_{j_{\iota}} \mid 1 \leq k \leq n,\ b_{j_{\iota}} \in \leadmon \},
\]
and we define 
\begin{equation} \label{eqn:QMonDef}
    Q =\ \partial \leadmon\, \cup\, \partial \trailmon .
\end{equation}

\begin{prop} \label{prop:effTestingProp}
     Let $I$ $=$ $(\mathcal{B})$ be a \SLI\ ideal. Then $I$ is efficient (that is, $I$ $=$ $(G)$) if and only if the following conditions hold: 
\begin{enumerate}
  \item[(i)] Every non-leading boundary monomial $b_j$ (\ie, $b_j$ $\in$ $\partial \mathcal{O} \setminus \leadmon$) is a multiple of at least one monomial in $Q$, and  
  \item[(ii)] $Q$ $\subseteq$ $(G)$.
\end{enumerate} 
\end{prop}

\proof
\noindent
$\mathbf{\Rightarrow}$: If $I$ $=$ $(G)$ is efficient, then $(G)$ contains each non-leading boundary monomial $b_j$; that is, $b_j$ $=$ $\sum_{b_{j_{\iota}}\in \leadmon} f_{{j_{\iota}}}\cdot g_{j_{\iota}}$, where the coefficients $f_{{j_{\iota}}}$ $\in$ $\gf[\mathbf{x}]$.  Recalling the form of the ideal generators $g_{j_{\iota}}$, it follows at once that the monomial $b_{j}$ is equal to a multiple of a monomial in $Q$; that is, (i) holds.  Furthermore, $\partial \trailmon$ consists of non-leading boundary monomials, so $\partial \trailmon$ $\subseteq$ $(G)$. Recalling (\ref{eqn:fDef}), we now see that  $x_k \cdot g_{j_{\iota}}$ $\in$ $(G)$ implies that each $x_k b_{j_{\iota}}$ $\in$ $(G)$; whence, $\partial \leadmon$ $\subseteq$ $(G)$, so $Q$ $\subseteq$ $(G)$; that is, (ii) holds.

\medskip

\noindent
$\mathbf{\Leftarrow}$: If conditions (i) and (ii) hold, we obtain at once that every non-leading boundary monomial is a member of $(G)$, so $I$ $=$ $(\mathcal{B})$ $=$ $(G)$ is efficient.  
\qed
\medskip

Given a \SLI\ ideal $I$, we can test it for efficiency by checking conditions (i) and (ii).  Condition (i) is straightforward, if possibly tedious, to check; it just depends on the order ideal and the choice of sets $\leadmon$ and $\trailmon$, which determine $Q$.  Condition (ii) can be tested as follows: One computes the $n \cdot \lambda$ products $x_k \cdot g_{j_{\iota}}$ ($g_{j_{\iota}}$ $\in$ $G$), and observes that the monomials appearing (non-trivially) in these products all lie in $Q$. Letting $E'$ be the set of indeterminates $\{ e_{\alpha,j_{\iota}} \mid 1\leq \alpha \leq n,\ b_{j_{\iota}} \in \leadmon \}$, we obtain a linear map 
\begin{equation} \label{eqn:mapSigma}
    \fn{\vartheta}{\Span_{\gf}(E') = \gf^{n \lambda}}{\Span_{\gf}(Q)}\ \text{given by } \ e_{k,j_{\iota}} \mapsto x_k \cdot g_{j_{\iota}},
\end{equation}
and condition (ii) holds if this map is surjective.  Accordingly, we say that $I$ is \textbf{$\vartheta$-efficient} whenever (i) holds and $\vartheta$ is surjective.  

\begin{rem} \label{rem:ThEffAndEff}
     Examples \ref{ssec:(1,6,6,10)} and \ref{(1,6,10,10,5)CaseC} exhibit efficient \SLI\ ideals $I$ that are not $\vartheta$-efficient, so $\vartheta$-efficiency is sufficient but not necessary for efficiency.  By contrast, the efficient ideal $(\mathcal{B})$ in Section \ref{ssec:RunningExampleII} is in fact $\vartheta$-efficient.  In this example, the domain of the linear map $\vartheta$ has dimension $n\cdot \lambda = 3\cdot 4 = 12$ and the codomain has dimension $10$, since the set $Q$ consists of the $10$ monomials of degree $3$ in $x_1,\, x_2,\, x_3$.  The elements of $G$ are then sufficiently general for $\vartheta$ to be surjective.
\end{rem}

Since the entries of the matrix of $\vartheta$ are the coefficients of the $g_{j_{\iota}}$, $\vartheta$-efficiency is an open condition on the \SLI\ locus $X_{\specl}$. That is, we have

\begin{cor}
    If $I$ $=$ $(\mathcal{B})$ $=$ $(G)$ is a $\vartheta$-efficient \SLI\ ideal, then there is an open set $[I]$ $\in$ $\mathcal{U}$ $\subseteq$ $X_{\specl}$ such that $[I']$ $\in$ $\mathcal{U}$ $\Rightarrow$ $I'$ is $\vartheta$-efficient. \qed
\end{cor}

\section{The tangent space at a point $[I]$ on the Hilbert Scheme} \label{sec:tanSpace}

\subsection{Tangent vectors at a point $[I]$ $\in$ ${\rm Hilb}^{\mu}_{\mathbb{A}^n_{\gf}}$} \label{ssec:tanSpGen}

	Let $I$ $\subseteq$ $\gf[\mathbf{x}]$ $=$ $R$ be an ideal of finite colength $\mu$ $=$ $\dim_{\gf}(R/I)$, and $\mathcal{B}$ $=$ $\{ g_1,\, \dots,\, g_{\nu}\}$ a border basis of $I$ with respect to an order ideal $\mathcal{O}$ $=$ $\{ t_1,\, \dots,\, t_{\mu} \}$.  
It is well-known that the tangent space $\tansp_{[I]}$ at the corresponding point $[I]$ $\in$ $\Hilb^{\mu}_{\mathbb{A}^n_{\gf}}$ is isomorphic to $\Hom{R}{I}{R/I}$ (see, \eg, \cite[Cor.\ 2.5, p.\ 13]{HartDefTh}).  
Hence, a tangent vector $\fn{v}{I}{R/I}$ at $[I]$ can be viewed as the vertical arrow in the following commutative diagram in which the top row is exact, $\phi(e_j)$ $=$ $g_j$, $1 \leq j \leq \nu$, and $\operatorname{Syz}$ $=$ the first syzygy module of the polynomials $ g_j $.  
\[
  \begin{diagram}
      {\rm Syz}  &\rInto^{i} & \oplus_{j=1}^{\nu} R\, e_j &\rOnto^{\phi} & I &\rTo & 0\\
         & & & \rdTo^{v'} &   \dTo^{v} & {}\\
                & & & &     R/I 
\end{diagram}
\]
It is therefore clear that a tangent vector $v$ corresponds to a choice of $\nu$ elements 
\[
    v'(e_j) \in R/I = \Span_{\gf}(\mathcal{O})
\]
such that for every tuple $(f_j)$ $\in$ ${\rm Syz}$ (and viewing the $v'(e_j)$ as elements of $R$), one has that 
\begin{equation} \label{eqn:tanVecCond}
    v'\left( \sum_{j=1}^{\nu}(f_j\, e_j)\right) = \sum_{j=1}^{\nu}f_j\, v'(e_j) \equiv 0 \mod{I};   
\end{equation}
moreover, it suffices for this condition to hold for every $(f_j)$ in a set of $R$-generators of ${\rm Syz}$.

Writing $v'(e_j)$ $=$ $\sum_{i=1}^{\mu} a_{ij}\, t_i$, we see that the tangent vector $v$ can be encoded as a $(\mu \nu)$-tuple of elements of $\gf$, as follows:
\begin{equation} \label{eqn:TanVec}
    v \leftrightarrow (a_{1,1},\, a_{2,1},\, \dots,\, a_{\mu,1},\, a_{1,2},\, a_{2,2}\, \dots,\, a_{\mu,2},\, a_{1,3},\, \dots,\, a_{\mu,\nu}). 
\end{equation}

A moment's reflection shows that, given any $f$ $\in$ $R$ and $\sum_{i=1}^{\mu}c_{i}t_i$ $\in$ $\Span_{\gf}(\mathcal{O})$, the product $f \cdot \left( \sum_{i=1}^{\mu}c_{i}t_i \right) $ reduces modulo $I$ to $\sum_{i=1}^{\mu} \mathbf{b}^f_i t_i$, where each coefficient $\mathbf{b}^f_i$ is a (unique) $\gf$-linear combination of the coefficients $c_i$.  Consequently, the sum $\sum_{j=1}^{\nu}f_j\, v'(e_j)$ in (\ref{eqn:tanVecCond}) reduces modulo $I$ to $\sum_{i=1}^{\mu}\mathbf{b}^{(f_j)}_{i}t_i$, where each of the coefficients $\mathbf{b}^{(f_j)}_i$ is a $\gf$-linear combination of the coefficients $a_{ij}$ that must vanish.  In other words, every syzygy $(f_j)$ imposes $\mu$ linear relations on the entries of the tuple (\ref{eqn:TanVec}) that must hold if the tuple is to encode a tangent vector; we will call these the \textbf{tangent space relations} associated to $(f_j)$.  As noted earlier, it suffices to check these conditions for each member of a set of $R$-generators of ${\rm Syz}$. 

Recalling from Section \ref{ssec:AlgorithmForNbrSyz} that a $\gf$-basis $\mathcal{L}$ of the linear syzygies provides a set of $R$-generators of ${\rm Syz}$, one sees from the foregoing that $\tansp_{[I]}$ is isomorphic to the $\gf$-vector subspace of tuples $(a_{ij})$ $\in$ $\gf^{\mu \nu}$ that satisfy all of the tangent space relations corresponding to the members of $\mathcal{L}$.  It is straightforward to compute these relations for specific examples via computer algebra; consequently, we can compute
\begin{equation} \label{eqn:tanSpDimComp}
    \dim_{\gf}(\tansp_{[I]}) = \mu \nu - \dim_{\gf}(\Span_{\gf}\{ \mathbf{b}^{(f_j)}_{i} \mid 1\leq i\leq \mu,\ (f_j) \in \mathcal{L} \}) .
\end{equation}


\subsection{Tangent vectors as $\gf[\epsilon]$-points} \label{ssec:kEpPts}

We write $\gf[\epsilon]$ for the dual numbers, that is, the $\gf$-algebra with $\epsilon^2$ $=$ $0$.  Recall that a tangent vector $v$ at $[I]$ can be viewed as a map of schemes $\fn{\theta_v}{\Spec(\gf[\epsilon])}{\Hilb^{\mu}_{\mathbb{A}^n_{\gf}}}$ such that the composition $\fnai{\Spec(\gf)}{\fnl{\theta_v}{\Spec(\gf[\epsilon])}{\Hilb^{\mu}_{\mathbb{A}^n_{\gf}}}}$ is the inclusion of the $\gf$-point $[I]$.  By the universal property of the Hilbert scheme, the map $\theta_v$ corresponds to a closed subscheme $Z_v$ $\subseteq$ $\Spec(\gf[\epsilon][\mathbf{x}])$ $=$ $R[\epsilon]$ such that $Z_v$ is finite and flat of degree $\mu$ over $\Spec(\gf[\epsilon])$ and the closed fiber is the closed subscheme $Z$ $\subseteq$ $\Spec(R)$ cut out by $I$.  The connection between this view of $v$ and the preceding, in which $v$ $\in$ $\Hom{R}{I}{O/I}$, is made as follows (see, \eg, \cite[prop. 2.3, p.\ 12]{HartDefTh}): the ideal $I_v$ $\subseteq$ $R[\epsilon]$ defining $Z_v$ has the form
\begin{equation} \label{eqn:tanVecForm}
	I_v = \{ f + \epsilon g \mid f \in I,\ g \in R, \text{ and } g \equiv v(f) \text{ mod } I \}.
\end{equation}


\section{An irreducible locus containing $[I]$} \label{sec:idealFam}

       Let $I$ be a \SLI\ ideal as in section \ref{sec:SLI}, from which we retain all notation.  We proceed to construct a map 
\[
    \fn{\EuScript{F}}{U}{ \Hilb^{\mu}_{\mathbb{A}^n_{\gf}}},
\]
where $U$ is an affine space and the image of $\EuScript{F}$ contains $[I]$ (indeed, $\EuScript{F}(U)$ contains the entire \SLI\ locus $X_{\specl}$).  In Section \ref{sec:derivMap} we compute the images of the standard unit tangent vectors at a $\gf$-point $p$ $\in$ $U$ under the derivative map $\fn{\EuScript{F}'_p}{\tansp_p}{\tansp_{[I_p]}}$, where $\EuScript{F}(p)$ $=$ $[I_{p}]$.  This will enable us to obtain lower bounds for the dimension of the image $\EuScript{F}(U)$, as explained in Section \ref{ssec:F(U)DimEstLemma}.

Roughly speaking, $\EuScript{F}(U)$ is obtained by ``translating'' $X_{\specl}$ around in $\Hilb^{\mu}_{\mathbb{A}^n_{\gf}}$ under maps $\fna{\Hilb^{\mu}_{\mathbb{A}^n_{\gf}}}{\Hilb^{\mu}_{\mathbb{A}^n_{\gf}}}$ induced by a family of automorphisms of $\mathbb{A}^{n}_{\gf}$, including the usual translations.


\subsection{Automorphisms of affine space} \label{ssec:isoms}

Let $A$ be a commutative and unitary ring, let $R_A$ $=$ $A[x_1,\dots,x_n]$ $=$ $A[\mathbf{x}]$, and write $\mathbb{A}^n_{A}$ $=$ $\Spec(R_A)$.  We define a family of automorphisms $\fna{\mathbb{A}^n_{A}}{\mathbb{A}^n_{A}}$ as follows: For $1$ $\leq$ $\alpha$ $\leq$ $n$, let $\Delta_{\alpha}$ be a finite set of monomials (including $1$) in the variables $\{ x_1, \dots, x_n \} \setminus \{ x_{\alpha} \}$ (specific choices of the sets $\Delta_{\alpha}$ are discussed in Sections \ref{ssec:ZDiscussion} and \ref{ssec:ZExtensionRem}).  We index each of these sets in some way, writing $\Delta_{\alpha}$ $=$ $\{ m_{\alpha,\delta} \mid 1 \leq \delta \leq |\Delta_{\alpha}| \}$. Note that the monomial $1$ will always have index $1$; that is, $m_{\alpha,1}$ $=$ $1$ for all $\alpha$.  For each choice of variable $x_{\alpha}$, monomial $m_{\alpha,\delta}$ $\in$ $\Delta_{\alpha}$, and scalar $z$ $=$ $z_{\alpha,\delta}$ $\in$ $A$, we obtain an automorphism
\[
     \fn{T^*_{ z_{\alpha,\delta}}}{R_A}{R_A},\ x_{\alpha} \mapsto x_{\alpha} + z_{\alpha,\delta}\cdot m_{\alpha,\delta} ,\ x_{\beta} \mapsto x_{\beta},\ \beta \neq \alpha.
\]
The map $T^*_{z_{\alpha,\delta}}$ induces an auomorphism $\fn{T_{ z_{\alpha,\delta}}}{\mathbb{A}^n_{A}}{\mathbb{A}^n_{A}}$.  For $m_{\alpha,1}$ $=$ $1$, the map $T_{z_{\alpha,\delta}}$ is just a translation of $\mathbb{A}^n_{A}$ in the $x_{\alpha}$-direction.

The ``translation'' of a subscheme $W$ $\subseteq$ $\mathbb{A}^n_{A}$ under an automorphism $\fn{T}{\mathbb{A}^n_{A}}{\mathbb{A}^n_{A}}$ is its pullback, denoted $W_T$: if $W$ $=$ $\Spec(R_A/\mathcal{I})$, then $W_T$ $=$ $\Spec(R_A/T^*(\mathcal{I}))$.   The following result is clear:

\begin{lem} \label{lem:isomLem}
     If the quotient $R_A/\mathcal{I}$ is generated as an $A$-module by a finite set $J$ $=$ $\{ f_1, \dots, f_d \}$ $\subseteq$ $R_A$ (resp.\ is $A$-free with basis $J$), then $R_A/T^*(\mathcal{I})$ is generated as an $A$-module by the set $ T^*(J)$ $=$ $\{ T^*(f_1) , \dots,  T^*(f_d) \}$ (resp. is $A$-free with basis $T^*(J)$).  \qed
\end{lem}


We list the elements of $\cup_{\alpha=1}^n \Delta_{\alpha}$ in the tuple 
\begin{equation} \label{eqn:Mdef}
    (m_{\alpha,\delta}) = \left( m_{1,1},\, m_{1,2},\, \dots,\, m_{1,|\Delta_1|},\,m_{2,1},\, m_{2,2}\, \dots,\, m_{n,|\Delta_n|} \right),
\end{equation}
and let $\mathbf{z}$ $=$ $(z_{\alpha,\delta})$ be a tuple of scalars  corresponding to the monomials in $\mathcal{M}$ in the order shown.  By composing the auomorphisms $T^*_{z_{\alpha,\delta}}$, we obtain the auomorphism of rings
\begin{equation}\label{eqn:f*_zDef}
    \fn{T^*_{\mathbf{z}} =  T^*_{z_{1,1}} \circ T^*_{z_{1,2}} \circ \dots \circ T^*_{z_{n,|\Delta_n|-1}} \circ T^*_{z_{n,|\Delta_n|}}}{R_A}{R_A},
\end{equation}
which in turn induces the automorphism of schemes
\begin{equation} \label{f_zDef}
    \fn{T_{\mathbf{z}}}{\mathbb{A}^n_{A}}{\mathbb{A}^n_{A}}.
\end{equation}


\subsection{Construction of the map $\EuScript{F}$} \label{ssec:Cconstr}

Let 
\[
    I = (\mathcal{B}) = (g_1, \dots, g_{\nu}) \subseteq \gf[x_1, \dots, x_n]
\]
be a \SLI\ ideal as in Section \ref{sec:SLI}.  In particular, the elements of the border basis $\mathcal{B}$ can be written as
\begin{equation} \label{eqn:gensOfStartingIdealI}
      g_j =                                \left\{ \begin{array}{l}
                            b_j, \text{ if } b_j \notin \leadmon,\vspace{.05in}\\
                            b_j - \sum_{C_{ij} \in \specl} c_{ij}t_i,\ c_{ij} \in \gf, \text{ if } b_j \in \leadmon \text{ ($\specl$  as in (\ref{eqn:SpecialC}))}
                                                \end{array} \right. .
\end{equation}

We introduce the set of variables 
\begin{equation} \label{eqn:V2def}
    \mathcal{Z} = \{ Z_{\alpha,\delta} \mid 1 \leq \alpha \leq n,\ 1 \leq \delta \leq |\Delta_{\alpha}| \}
\end{equation}
corresponding to the scalars $z_{\alpha,\delta}$ introduced in Section \ref{ssec:isoms}.  Let $A$ denote the polynomial ring $\gf[\specl, \mathcal{Z}]$, and let the ideal $\mathfrak{I}$ $\subseteq$ $A[\mathbf{x}]$ be generated by the $\mathcal{O}$-border prebasis
\begin{equation}\label{eqn:preBasisOfGs}
    \mathcal{B}_{\mathfrak{I}} = \left\{ G_1, \dots, G_{\nu} \right\}, \ G_j = \left\{ \begin{array}{l}
                            b_j, \text{ if } b_j \notin \leadmon,\vspace{.05in}\\
                            b_j - \sum_{C_{ij} \in \specl} C_{ij}\, t_i,\ \text{ if } b_j \in \leadmon
                                                \end{array} \right..
\end{equation}

\begin{lem} \label{lem:preBasisIsBasis}
    The prebasis $\mathcal{B}_{\mathfrak{I}}$ is in fact an $\mathcal{O}$-border basis, and the quotient $A[\mathbf{x}]/\mathfrak{I}$ is $A$-free with basis $\mathcal{O}$.
\end{lem}

\proof         
Arguing as in the proof of Proposition \ref{prop:specialCCor}, one sees that for all pairs of neighbors $b_j$ and $b_{j'}$ $\in$ $\partial \mathcal{O}$, the polynomial $S(G_j,G_{j'})$ of Section \ref{ssec:NbrSyz} reduces to $0$ modulo $\mathcal{B}_{\mathfrak{I}}$.  Proposition \ref{prop:NbrSyzFacts} now implies that $\mathcal{B}_{\mathfrak{I}}$ is an $\mathcal{O}$-border basis, so $A[\mathbf{x}]/\mathfrak{I}$ is $A$-free with basis $\mathcal{O}$.
\qed
\medskip

Lemmas \ref{lem:isomLem} and \ref{lem:preBasisIsBasis} yield that the quotient $A[\mathbf{x}]/T^*_{\mathbf{Z}}(\mathfrak{I})$ is $A$-free with basis $T^*_{(\mathbf{Z})}(\mathcal{O})$; consequently, the induced map
\begin{equation} \label{eqn:univFamOverU}
    \fna{\Spec(A[\mathbf{x}]/T^*_{\mathbf{Z}}(\mathfrak{I}))}{\Spec(A) = U}
\end{equation}
is finite of degree $\mu$ $=$ $|\mathcal{O}|$ and flat, and so, by the universal property of the Hilbert scheme, corresponds to a map $\fn{\EuScript{F}}{U}{ \Hilb^{\mu}_{\mathbb{A}^n_{\gf}}}$.  

Let $p$ $=$ $(\mathbf{c},\mathbf{z})$ be a $\gf$-point of $U$ $=$ $\Spec(\gf[\specl,\mathcal{Z}])$, where $\mathbf{c}$ $=$ $(c_{ij})$ is a tuple of elements of $\gf$ indexed by the \SLI\ index-pairs $(i,j)$, and $\mathbf{z}$ $=$ $(z_{\alpha,\delta})$ is a tuple of elements of $\gf$ corresponding to the monomials $m_{\alpha,\delta}$ and ordered as in (\ref{eqn:Mdef}).  Let $\fn{T^*_{\mathbf{z}}}{\gf[\mathbf{x}]}{\gf[\mathbf{x}]}$ be the corresponding map (\ref{f_zDef}).  A moment's reflection shows that the fiber of (\ref{eqn:univFamOverU}) over $p$ is the closed subscheme $\gf[\mathbf{x}]/T^*_{\mathbf{z}}(I_{\mathbf{c}})$, where $I_{\mathbf{c}}$ is the \SLI\ ideal
\[
    I_{\mathbf{c}} = (g'_1, g'_2, \dots, g'_{\nu}),\  g'_j = \left\{ \begin{array}{l}
                            b_j, \text{ if } b_j \notin \leadmon,\vspace{.05in}\\
                            b_j - \sum_{C_{ij} \in \specl} c_{ij}\, t_i, \text{ if } b_j \in \leadmon
                                                                   \end{array} \right. .
\]
In particular, it is clear that the fiber over the origin $(\mathbf{0},\mathbf{0})$ $\in$ $U$ is $\Spec(\gf[\mathbf{x}]/I_0)$, where $I_0$ is (as defined following Corollary \ref{cor:DistIdealsSuppAtOnePt}) the monomial ideal $(\partial \mathcal{O})$, and the fiber over the point $((c_{ij}),\mathbf{0})$ is $\Spec(\gf[\mathbf{x}]/I_{\mathbf{c}})$, so that $\EuScript{F}((\mathbf{0},\mathbf{0}))$ $=$ $[I_0]$ and $\EuScript{F}(((c_{ij}),\mathbf{0}))$ $=$ $[I_{\mathbf{c}}]$.

\begin{rem} \label{rem:supportIsPtRem}
    Since the ideal $I_{\mathbf{c}}$ is supported at one point (the origin) of $\mathbb{A}^n_{\gf}$ by Corollary \ref{cor:DistIdealsSuppAtOnePt}, it follows that $T^*_{\mathbf{z}}(I_{\mathbf{c}})$ is also supported at one point of $\mathbb{A}^n_{\gf}$.  Consequently, every point $[I']$ $\in$ $\EuScript{F}(U)$ corresponds to an ideal $I'$ $\subseteq$ $\gf[\mathbf{x}]$ that is supported at one point.
\end{rem}


\subsection{Finding lower bounds for $\dim(\EuScript{F}(U))$} \label{ssec:F(U)DimEstLemma}

Our method for bounding the dimension of $\EuScript{F}(U)$ from below is summarized by the following elementary

\begin{prop} \label{prop:F(U)DimLowerBd}

Let $p$ be a $\gf$-point of $U$.  Suppose given a set of tangent vectors $\{ v_1, \dots, v_L \}$ $\subseteq$ $\tansp_p$ such that the image set 
\[
    \{ \EuScript{F}'_p(v_1),\dots, \EuScript{F}'_p(v_L) \} \subseteq \tansp_{[I_p]}
\]
is linearly independent.  Then $L$ $\leq$ $\dim(\EuScript{F}(U))$.

\end{prop} 

\proof
        Let $\dim(U)$ $=$ $\dim_{\gf}(\tansp_p)$ $=$ $d$.  The hypothesis implies that $\dim_{\gf}(\EuScript{F}'_p(\tansp_p))$ $\geq$ $L$; whence, $\dim(\ker(\EuScript{F}'_p))$ $\leq$ $d-L$.  It follows that any component of the fiber $\EuScript{F}^{-1}([I_p])$ through $p$ has dimension $\leq$ $d-L$, so by the theorem on the dimension of fibers of a morphism, $d-L$ $\geq$ $d - \dim(\EuScript{F}(U))$ $\Rightarrow$ $L$ $\leq$ $\dim(\EuScript{F}(U))$, as asserted.
\qed
\medskip

Recall that an irreducible component of $\Hilb^{\mu}_{\mathbb{A}^n_{\gf}}$ is called \textbf{elementary} if every point $[I']$ on it parameterizes a subscheme $\Spec(\gf[\mathbf{x}]/I')$ that is concentrated at one point. By Remark \ref{rem:supportIsPtRem}, this property holds for all $[I']$ $\in$ $\EuScript{F}(U)$. Also recall that if $[I']$ is a smooth point on an elementary component, then we call $I'$ a \textbf{generic} ideal.  
Proposition \ref{prop:F(U)DimLowerBd} leads to the following simple criterion for identifying elementary components and generic ideals:

\begin{prop} \label{prop:genericCond}
Let $L$ be a lower bound for $\dim(\EuScript{F}(U))$ as in Proposition {\rm \ref{prop:F(U)DimLowerBd}}. If there is a point $[I]$ $\in$ $\EuScript{F}(U)$ such that 
$\dim_{\gf}(\tansp_{[I]})$ $=$ $L$, then the closure $\overline{\EuScript{F}(U)}$ is an elementary component of $\Hilb^{\mu}_{\mathbb{A}^n_{\gf}}$ of dimension $L$ on which $[I]$ is a smooth point; consequently, $I$ is a generic ideal.
\end{prop}

\proof
     The hypothesis implies that $L$ is both a lower bound and an upper bound for $\dim(\EuScript{F}(U))$; whence, $\dim(\overline{\EuScript{F}(U)})$ $=$ $L$ and $[I]$ is a smooth point on $\overline{\EuScript{F}(U)}$.  Then the unique irreducible component of $\Hilb^{\mu}_{\mathbb{A}^n_{\gf}}$ containing $[I]$ must be $\overline{\EuScript{F}(U)}$, which is accordingly an elementary component, and we are done.
\qed
\medskip


\section{The derivative map $\fn{\EuScript{F}'_p}{\tansp_{p}}{\tansp_{[I_p]}}$} \label{sec:derivMap}
Let $p$ $=$ $(\mathbf{c},\mathbf{z})$ $\in$ $U$. In this section of the paper, we study the derivative map $\fn{\EuScript{F}'_p}{\tansp_p}{\tansp_{[I_p]}}$.  Since $U$ $=$ $\Spec(\gf[\specl,\mathcal{Z}])$, a basis of $\tansp_p$ is given by unit vectors in the directions corresponding to the indeterminates $C_{ij}$ $\in$ $\specl$ and $Z_{(\alpha,\delta)}$ $\in$ $\mathcal{Z}$.  Let $X$ denote one of these variables, let $X'$ $\neq$ $X$ stand for any of the others, and let $p_X$, $p_{X'}$ $\in$ $\gf$ denote the corresponding components of $p$.  Then a unit vector in the $X$-direction at $p$ is given by the map
\[
      \fn{v_{p,X}}{\Spec(\gf[\epsilon])}{U}\ \text{defined by } X \mapsto p_X + \epsilon,\ X' \mapsto p_{X'}.
\]
The image $\EuScript{F}'_p(v_{p,X})$ $\in$ $\tansp_{[I_p]}$ is then the map
\[
      \EuScript{F}'_p(v_{p,X}): \Spec(\gf[\epsilon]) \stackrel{v_{p,X}}{\rightarrow} U  \stackrel{\EuScript{F}}{\rightarrow} \Hilb^{\mu}_{\mathbb{A}^n_{\gf}},
\]
which corresponds to an ideal $I_{p,X}$ $\subseteq$ $\gf[\epsilon][\mathbf{x}]$ such that the quotient $\gf[\epsilon][\mathbf{x}]/I_{p,X}$ is $\gf[\epsilon]$-free of rank $\mu$.  Recall that $\EuScript{F}$ is defined by the ideal 
\[
     T^*_{\mathbf{Z}}(\mathfrak{I}) = (\{T^*_{\mathbf{Z}}(G_j) \mid 1 \leq j \leq \nu \}) \subseteq \gf[\specl,\mathcal{Z}][\mathbf{x}],
\]
where the $G_j$ are defined in equation (\ref{eqn:preBasisOfGs}). 
Thus, $I_{p,X}$ is the image of this ideal under the substitutions $X \mapsto p_X + \epsilon$, $X' \mapsto p_{X'}$.  

We now restrict attention to the point
\begin{equation} \label{eqn:pRestriction}\
     p = ((c_{ij}),\mathbf{0}),
\end{equation}
so that $I_p$ is \SLI\ with border basis  $\mathcal{B}$ as shown in (\ref{eqn:gensOfStartingIdealI}).
We proceed to evaluate the tangent vectors $\EuScript{F}'_p(v_{p,X})$ $\in$ $\tansp_{[I_p]}$ for each of the cases $X$ $=$ $C_{ij}$ $\in$ $\specl$ and $X$ $=$ $Z_{\alpha,\delta}$ $\in$ $\mathcal{Z}$.


\subsection{The tangent vectors $\EuScript{F}'_p(v_{p,C_{ij}})$ for $C_{ij}$ $\in$ $\specl$} \label{ssec:TVFam1}


If $X$ $=$ $C_{ij}$ $\in$ $\specl$, and $p$ is as in (\ref{eqn:pRestriction}), one sees easily that for all $j' \neq j$, $1$ $\leq$ $j'$ $\leq$ $\nu$, the image of $T^*_{\mathbf{Z}}(G_{j'})$ under $C_{ij}$ $\mapsto$ $p_{C_{ij}} + \epsilon$, $X' \mapsto p_{X'}$, is $g_{j'}$, and the image of $T^*_{\mathbf{Z}}(G_{j})$ is $g_{j} - \epsilon\, t_{i}$, where $g_j$ and $g_{j'}$ are as in (\ref{eqn:gensOfStartingIdealI}). According to (\ref{eqn:tanVecForm}), the tangent vector $\EuScript{F}'_p(v_{p,C_{ij}})$ $=$ $v_{p,ij}$ corresponds to the element of $\Hom{R}{I_p}{R/I_p}$ given by $g_{j}$ $\mapsto$ $(-t_{i}\text{ mod }I_p) = -t_{i}$, and $g_{j'}$ $\mapsto$ $0$ for $j'$ $\neq$ $j$.  The corresponding tuple $(a_{i'j'})$ (equation (\ref{eqn:TanVec})) has all components equal to $0$ except for $a_{ij}$ $=$ $-1$.  The following lemma is immediate:  

\begin{lem} \label{lem:YtanVecsLinInd}
     Let $p$ be a point as in equation {\rm (\ref{eqn:pRestriction})}. Then the family of tangent vectors 
\[
    \EuScript{S}_p = \{ v_{p,ij} \mid C_{ij} \in \specl \} \subseteq \tansp_{[I_p]}
\]
is $\gf$-linearly independent and of cardinality 
\[
    |\EuScript{S}_p| =  |\leadmon| \cdot |\trailmon| = \lambda \cdot \tau.\ 
\] \qed
\end{lem}

\begin{rem} \label{rem:spineRem}
    It is clear that $\EuScript{S}_p$ is a basis of the tangent space to the \SLI\ locus $X_\specl$ $\subseteq$ $\mathbb{B}_{\mathcal{O}}$ at the point $[I_p]$.
\end{rem}


\subsection{The tangent vectors $\EuScript{F}'_p(v_{p, Z_{\alpha,\delta}})$ for $Z_{\alpha,\delta}$ $\in$ $\mathcal{Z}$} \label{ssec:TVFam2}

Now consider the case $X$ $=$ $Z_{\alpha,\delta}$ for some $1$ $\leq$ $\alpha$ $\leq$ $n$, $m_{\alpha,\delta}$ $\in$ $\Delta_{\alpha}$ (recall that the latter is a finite set of monomials not involving $x_{\alpha}$).  Recalling that $p$ is a point as in equation {\rm (\ref{eqn:pRestriction})}, one sees that the ideal $I_{p,X}$ is obtained by applying the substitutions $X$ $\mapsto$ $p_X + \epsilon$ $=$ $\epsilon$, $X'$ $\mapsto$ $p_{X'}$, to the polynomials $T^*_{\mathbf{Z}}(G_j)$. Under these substitutions, which amount to replacing $x_{\alpha}$ by $x_{\alpha} + \epsilon \cdot m_{\alpha, \delta}$ in the polynomials $g_j$, one has that
\[
        T^*_{\mathbf{Z}}(G_j) \mapsto g_j + \epsilon \cdot \frac{\partial g_j}{\partial x_{\alpha}} \cdot m_{\alpha,\delta},\ 1 \leq j \leq \nu.
\]
Hence, (\ref{eqn:tanVecForm}) yields that the tangent vector $\EuScript{F}'_p(v_{p, Z_{\alpha,\delta}})$ $=$ $v_{p,\alpha, \delta}$ corresponds to the element of $\Hom{R}{I_p}{R/I_p}$ given by 
\begin{equation} \label{eqn:ZTanVecHomo}
     g_j\ \mapsto \frac{\partial g_j}{\partial x_{\alpha}} \cdot m_{\alpha,\delta}\ \equiv \ \sum_{i=1}^{\mu}a_{ij}t_i\ (\text{mod } I_p),\ \ 1 \leq j \leq \nu.
\end{equation}
For a point $p$ as in (\ref{eqn:pRestriction}), we let

\begin{equation} \label{eqn:tanVecSetZ}
    \EuScript{Z}_p = \{ v_{p,\alpha,\delta} \mid 1 \leq 
\alpha \leq n,\ 1 \leq \delta \leq |\Delta_{\alpha}| \} \subseteq \tansp_{I_p}.
\end{equation}


\section{Lex-segment complement order ideals} \label{sec:spCaseRestr}

In this section we discuss the order ideals and associated \SLI\ ideals that are used in all of our examples in Section \ref{sec:examples}.

\subsection{Definition} \label{ssec:lexSegComps}
From this point on, a monomial inequality (such as $m_1$ $>$ $m_2$) shall be with respect to the lexicographic order with $x_1$ $>$ $x_2$ $>$ \dots $>$ $x_n$.  Let $\EuScript{L}$ be a (proper) lex-segment ideal of finite colength in $R$ $=$ $\gf[\mathbf{x}]$ (see, \eg, \cite[Sec.\ 5.5.B, p.\ 258]{KreutzerAndRobbianoVolTwo}), and let $\mathcal{O}$ be the set of monomials that are not in $\EuScript{L}$; we call $\mathcal{O}$ a \textbf{lex-segment complement} order ideal.  Writing $R_d$ $\subseteq$ $R$, 
$\EuScript{L}_d$ $\subseteq$ $\EuScript{L}$, and $\mathcal{O}_d$ $\subseteq$ $\mathcal{O}$ for the subsets of monomials of degree $d$, we let $m_d$ denote the lex-minimum element of $\EuScript{L}_d$ when this set is non-empty, and $s$ $\geq$ $r$ $>$ $0$ the integers such that 
\[
    \mathcal{O}_d \ =\ \left\{
    \begin{array}{l}
          R_d, \text{ if } 0 \leq d < r;\vspace{.05in}\\
         \{m \in R_d \mid  m < m_d \} \neq \emptyset, \text{ if } r\leq d \leq s;\vspace{.05in}\\
         \emptyset, \text{ if } d > s.
    \end{array}
                      \right. 
\]
Note that for $r$ $\leq$ $d$ $<$ $s$ one has that $m_{d+1}$ $\leq$ $x_n \cdot m_d$.  Here is a simple example in $3$ variables with Hilbert function $(1,3,2,1)$; the basis monomials are underlined and the boundary monomials are shown in boldface.   In this case, $r=2$, $s=3$, $m_2$ $=$ $x_2^2$, and $m_3$ $=$ $x_2 x_3^2$:
\[
    \begin{array}{c}
        \underline{1}\vspace{.05in}\\
        \underline{x_1}\ \ \underline{x_2}\ \ \underline{x_3}\vspace{.05in}\\
        \mathbf{x_1^2}\ \ \mathbf{x_1 x_2} \ \ \mathbf{x_1 x_3}\ \ \mathbf{x_2^2}\ \ \underline{x_2 x_3}\ \ \underline{x_3^2}\vspace{.05in}\\
        x_1^3\ \ x_1^2 x_2\ \  x_1^2 x_3\ \ x_1 x_2^2\ \ \mathbf{x_1 x_2 x_3}\ \ \mathbf{x_1 x_3^2}\ \ x_2^3\ \ \mathbf{x_2^2 x_3}\ \ \mathbf{x_2 x_3^2}\ \ \underline{x_3^3}\vspace{.05in}\\
    x_1^4\ \ x_1^3 x_2\ \ \ \ \ \ \ \ \dots\ \ \ \ \ \ \ \ x_1x_2 x_3^2\ \ \mathbf{x_1 x_3^3}\ \  x_2^4\ \ \ \ \ \ \dots\ \ \ \ \ \ \ x_2^2 x_3^2 \ \mathbf{x_2 x_3^3} \ \mathbf{x_3^4}
   \end{array}
\]

\begin{lem} \label{lem:lexSegMonomLem}
    If $\mathcal{O}$ is a lex-segment complement order ideal, $m$ is a monomial $\notin \mathcal{O}$, and $m'$ is a monomial such that $d$ $=$ $\deg(m')$ $\geq$ $\deg(m)$ and $m'$ $>$ $m$, then $m'$ $\notin$ $\mathcal{O}$.
\end{lem}

\proof
    Let $\deg(m')-\deg(m)$ $=$ $u \geq 0$.  Then $m''$ $=$ $m \cdot x_n^u$ is the lex-minimum monomial of degree  $d$ that is $\geq m$; whence, $m'$ $\geq$ $m''$.  Furthermore, the hypothesis on $\mathcal{O}$ implies that $m''$ $\notin$ $\mathcal{O}$, which in turn yields $m'$ $\notin$ $\mathcal{O}$, as desired.
\qed

 
\subsection{The sets of monomials $\Delta'_{\alpha}$} \label{ssec:ZDiscussion}

Recall that in order to define the map $\fn{\EuScript{F}}{U}{\Hilb^{\mu}_{\mathbb{A}^n_{\gf}}}$, we must choose finite sets of monomials $\Delta_{\alpha}$ as in Section \ref{ssec:isoms}.  Here we describe the particular sets $\Delta'_{\alpha}$ that we most often use when $\mathcal{O}$ is a lex-segment complement. 

We claim that for every variable $x_{\alpha}$, there is a smallest exponent $e'_{\alpha}$ $\geq$ $0$ of $x_n$ such that $x_{\alpha}\cdot x_n^{e'_{\alpha}}$ $\notin$ $\mathcal{O}$ but $x_n^{e'_{\alpha}}$ $\in$ $\mathcal{O}$.  When $\alpha$ $=$ $n$, it is clear that $e'_n$ is the largest exponent $e$ such that $x_n^e$ $\in$ $\mathcal{O}$.  We then observe that for any $\alpha$ $\neq$ $n$, $x_{\alpha}\cdot x_n^{e'_n}$ $>$ $x_n\cdot x_n^{e'_n}$ $\Rightarrow$ $x_{\alpha} \cdot x_n^{e'_n}$ $\notin$ $\mathcal{O}$ by Lemma \ref{lem:lexSegMonomLem}, but $x_n^{e'_n}$ $\in$ $\mathcal{O}$; from this it follows that $e'_{\alpha}$ exists and is $\leq$ $e'_n$.   We define
\begin{equation} \label{eqn:spCaseBAlphaDefn}
    b_{j_{\alpha}} = x_{\alpha} \cdot x_n^{e'_{\alpha}} \in \partial \mathcal{O},\ \ t_{i_{\alpha}} = x_n^{e'_{\alpha}} \in \mathcal{O}.
\end{equation}
Note that $b_{j_{\alpha}}$ can be characterized as the lex-minimum boundary monomial that is divisible by $x_{\alpha}$.  

We now choose the monomial sets $\Delta'_{\alpha}$, $1$ $\leq$ $\alpha$ $\leq$ $n$, as follows:
\begin{equation} \label{eqn:DeltaDelta'RestrDef}
   \begin{array}{c}
            \Delta'_n = \{ m_{n,1}\} = \{1\}, \text{ and, for } 1 \leq \alpha \leq n-1,\vspace{.05in}\\
            \Delta'_{\alpha} = \left\{ m_{\alpha,\delta} \in \gf[x_{\alpha+1},\dots,x_n] 
            \left|  \begin{array}{l}
                        t_{i_{\alpha}} m_{\alpha,\delta} \in \mathcal{O} \setminus \trailmon, \text{ if } b_{j_{\alpha}} \in \leadmon,\vspace{.05in}\\
              t_{i_{\alpha}} m_{\alpha,\delta} \in \mathcal{O}, \text{ if } b_{j_{\alpha}} \notin \leadmon 
                    \end{array} 
             \right. 
            \right\}.
   \end{array}
\end{equation}
One checks easily that $1$ $=$ $m_{\alpha,1}$ $\in$ $\Delta'_{\alpha}$ for all $\alpha$.  Note that we have replaced the condition $\Delta'_{\alpha}$ $\subseteq$ $\gf[x_1, \dots, \widehat{x_{\alpha}}, \dots, x_n]$ of Section \ref{ssec:isoms} with the seemingly stricter $\Delta'_{\alpha} \subseteq \gf[x_{\alpha+1}, \dots, x_n]$, since if $m_{\alpha,\delta}$ is divisible by one of $x_1$, \dots, $x_{\alpha-1}$, then, by Lemma \ref{lem:lexSegMonomLem},
\[
    t_{i_{\alpha}} \cdot m_{\alpha,\delta} > t_{i_{\alpha}}\cdot x_{\alpha} = b_{j_{\alpha}} \notin{O}\ \Rightarrow\ t_{i_{\alpha}} \cdot m_{\alpha,\delta} \notin \mathcal{O}.
\]

As in (\ref{eqn:V2def}), we write 
\[
    \mathcal{Z}' = \{ Z_{\alpha,\delta} \mid 1 \leq \alpha \leq n,\ 1 \leq \delta \leq |\Delta'_{\alpha}| \},
\]
and for $p$ as in (\ref{eqn:pRestriction}), we denote the set of tangent vectors (\ref{eqn:tanVecSetZ}) associated to  $\mathcal{Z}'$ by

\[
   \EuScript{Z}'_p\ =\  \{ v_{p,\alpha,\delta} \mid 1 \leq \alpha \leq n,\ m_{\alpha,\delta} \in \Delta'_{\alpha}  \}\ \subseteq\ \tansp_{[I_p]}.
\]


\subsection{Ideals $I_p$ for which the set $\EuScript{S}_p$ $\cup$ $\EuScript{Z}'_p$ $\subseteq$ $\tansp_{[I_p]}$ is linearly independent} 
\label{ssec:ZMainResult}

The key technical result enabling us to find elementary components of $\Hilb^{\mu}_{\mathbb{A}^n_{\gf}}$ is the following

\begin{prop} \label{prop:GeneralZProp}  
    Let $\mathcal{O}$ $\neq$ $\{ 1 \}$ be a lex-segment complement order ideal, and let $p$ be as in {\rm (\ref{eqn:pRestriction})}, so that $I_p$ is \SLI\ with $\mathcal{O}$-border basis $\mathcal{B}$ $=$ $\{ g_j \mid 1 \leq j \leq \nu \}$ as in { \rm (\ref{eqn:gensOfStartingIdealI})}. Suppose that $\mathcal{B}$ has the property that the boundary monomial $b_{j_{\alpha}}$ {\rm (\ref{eqn:spCaseBAlphaDefn})} is the lex-leading monomial in $g_{j_{\alpha}}$ for all $1$ $\leq$ $\alpha$ $\leq$ $n$.  Then the set $\EuScript{S}_p \cup \EuScript{Z'}_p$ $\subseteq$ $\tansp{[I_p]}$ is $\gf$-linearly independent.  
\end{prop}

\proof
    Suppose we have a relation
\begin{equation} \label{eqn:Z2LinComb}
    \sum_{v_{p,ij} \in \EuScript{S}_p} d_{ij}\,v_{p,ij} + \sum_{v_{p,\alpha,\delta} \in \EuScript{Z'}_p} d_{\alpha,\delta}\, v_{p,\alpha,\delta}\ =\ 0,\ \  \ d_{ij},\ d_{\alpha,\delta} \in \gf . 
\end{equation} 
We must show that all the coefficients in the linear combination vanish.

By (\ref{eqn:ZTanVecHomo}), we have that the homomorphism $\fna{I_p}{R/I_p}$  corresponding to $v_{p,\alpha,\delta}$ $\in$ $\EuScript{Z'}_p$ is given by
\[
     g_j\ \mapsto \frac{\partial g_j}{\partial x_{\alpha}} \cdot m_{\alpha,\delta}\ (\text{mod } I_p),\ \ 1 \leq j \leq \nu.
\]

For each $1 \leq \alpha \leq n$ and each $m_{\alpha,\delta} \in \Delta'_{\alpha}$, let 
\[
    t_{i_{\alpha,\delta}} = t_{i_{\alpha}} \cdot m_{\alpha,\delta} \in \mathcal{O}.
\]
Our first goal is to show that the $(i_{\alpha,\delta},j_{\alpha})$-component of $v_{p,\alpha,\delta}$ is non-zero for all $1 \leq \alpha \leq n$ and all $1 \leq \delta \leq |\Delta'_{\alpha}|$.  This component is the coefficient of $t_{i_{\alpha,\delta}}$ in $\frac{\partial g_{j_{\alpha}}}{\partial x_{\alpha}} \cdot m_{\alpha,\delta}$ mod $I_p$. Consider first the case in which $b_{j_{\alpha}}$ $\notin$ $\leadmon$, so that $g_{j_{\alpha}}$ $=$ $b_{j_{\alpha}}$ $=$ $x_{\alpha}x_n^{e'_{\alpha}}$.  In this case, 
\[
    \frac{\partial g_{j_{\alpha}}}{\partial x_{\alpha}} = \frac{\partial b_{j_{\alpha}}}{\partial x_{\alpha}} = \left\{
        \begin{array}{l}
            x_n^{e'_{\alpha}} = t_{i_{\alpha}},\ \text{if } 1 \leq \alpha < n,\vspace{.05in}\\
            (e'_{n}+1) x_n^{e'_n} = (e'_n+1) \cdot t_{i_{\alpha}},\ \text{if } \alpha = n.
        \end{array} \right.
\]
It follows that 
\[
    \frac{\partial g_{j_{\alpha}}}{\partial x_{\alpha}} \cdot m_{\alpha, \delta} = (1 \text{ or } e'_{\alpha}+1) \cdot t_{i_{\alpha,\delta}} \equiv (1 \text{ or } e'_{\alpha}+1) \cdot t_{i_{\alpha,\delta}} \mod I_p, 
\]
so in either case the $(i_{\alpha,\delta},j_{\alpha})$-component of $v_{p,\alpha,\delta}$ (either $1$ or $e'_{\alpha}+1$) is non-zero, since $\operatorname{char}(\gf)$ $=$ $0$.

Next consider the case in which $b_{j_{\alpha}}$ $\in$ $\leadmon$, so that $g_{j_{\alpha}}$ $=$ $b_{j_{\alpha}}-N_{j_{\alpha}}$, where $N_{j_{\alpha}}$ $\in$ $\Span_{\gf}(\trailmon)$.  As before, $\frac{\partial b_{j_{\alpha}}}{\partial x_{\alpha}}$ $=$ $(1 \text{ or } e'_{n}+1) \cdot t_{i_{\alpha}}$, which contributes a non-zero multiple of $t_{i_{\alpha,\delta}}$ to $(\frac{\partial g_{j_{\alpha}}}{\partial x_{\alpha}} \cdot m_{\alpha,\delta}$ mod $I_p)$. In fact, this is the only (non-zero) contribution to the $(i,j_{\alpha})$-component of $v_{p,\alpha,\delta}$ for any $i$, because we have the following claim: 
\begin{equation}\label{txt:keyClaim}
  \begin{array}{l}
    \text{\emph{If $b_{j_{\alpha}}$ $\in$ $\leadmon$ and $m$ is a monomial appearing non-trivially}}\\
    \text{\emph{in $N_{j_{\alpha}}$, then $m$ is not divisible by $x_{\alpha}$.}}
  \end{array}
\end{equation}

To prove the claim, suppose that $x_{\alpha}$ divides $m$.  By hypothesis, we have that $b_{j_{\alpha}}$ $=$ $x_{\alpha} x_n^{e'_{\alpha}}$ $>$ $m$, which implies that $m$ $=$ $x_{\alpha}\cdot x_n^e$, with $e'_{\alpha}$ $>$ $e$.  But $m$ $\in$ $\trailmon$, which implies that $m \cdot x_{n}$ $\in$ $\partial \mathcal{O}$ $\setminus$ $\leadmon$.  Consequently, we either have that $b_{j_{\alpha}}$ $=$ $m\cdot x_n$, which contradicts $b_{j_{\alpha}}$ $\in$ $\leadmon$, or $b_{j_{\alpha}}$ $=$ $x_n^w\cdot (m\cdot x_n)$ with $w \geq 1$, which contradicts that the leading monomial $b_{j_{\alpha}}$ is a minimal boundary monomial.  We conclude that $m$ cannot be divisible by $x_{\alpha}$, as claimed. 

 So far we have established that the $(i_{\alpha,\delta},j_{\alpha})$-component of the tangent vector $v_{p,\alpha,\delta}$ is non-zero for all $\alpha$ and all $m_{\alpha,\delta}$ $\in$ $\Delta'_{\alpha}$.
We now show by descending induction on $\alpha$ that the coefficients $d_{\alpha,\delta}$ in (\ref{eqn:Z2LinComb}) are all equal to $0$.  

We begin with the case $\alpha$ $=$ $n$ (recall that $\Delta'_n=\{ 1\}$).  We claim that for all tuples $(\beta,\delta')$ $\neq$ $(n,1)$, which implies that $\beta$ $<$ $n$, the $(i_{n,1},j_n)$-component of $v_{p,\beta,\delta'}$ is $0$.  This component is the coefficient of 
\[
    t_{i_{n,1}} = t_{i_n}\cdot 1 = t_{i_n} = x_n^{e'_n} \text{ in }
 \left( \frac{\partial g_{j_n}}{\partial x_{\beta}} \cdot m_{\beta,\delta'} \mod I_p \right).
\]
In case $b_{j_n}$ $\notin$ $\leadmon$, then $g_{j_n}$ $=$ $b_{j_n}$ $=$ $x_n^{e'_n+1}$, so $\frac{\partial g_{j_n}}{\partial x_{\beta}}$ $=$ $0$. In case $b_{j_n}$ $\in$ $\leadmon$, any monomial $m$ appearing non-trivially in $N_{j_{n}}$ must satisfy $x_n^{e'_n+1}$ $>$ $m$ and that $m$ is not divisible by $x_n$, by (\ref{txt:keyClaim}); whence, $m = 1$ $\in$ $\trailmon$, so $x_{\alpha} \cdot 1$ $\in$ $\partial \mathcal{O}$ for all $\alpha$, and we are in the excluded case $\mathcal{O}$ $=$ $\{ 1 \}$.  It follows once again that $g_{j_n}$ $=$ $b_{j_n}$ $\Rightarrow$ $\frac{\partial g_{j_n}}{\partial x_{\beta}}$ $=$ $0$. We conclude that the $(i_{n,1},j_n)$-component of $v_{p,\beta,\delta'}$ is $0$ for all $\beta < n$ and $1 \leq \delta' \leq |\Delta'_{\beta}|$, as claimed.

We next note:
\begin{equation}\label{txt:nonSpecialIndices}
    \text{\emph{None of the index pairs $(i_{\alpha,\delta},j_{\alpha})$ are \SLI.}}
\end{equation}
To see this, recall from (\ref{eqn:SpecialC}) that the index pair $(i,j)$ is \SLI\ if and only if $b_j$ $\in$ $\leadmon$ and $t_i$ $\in$ $\trailmon$.  However, by the definition (\ref{eqn:DeltaDelta'RestrDef}) of our sets $\Delta'_{\alpha}$, we have that $t_{i_{\alpha,\delta}}$ $\notin$ $\trailmon$ whenever $b_{j_{\alpha}}$ $\in$ $\leadmon$. 

The foregoing implies that in the sum (\ref{eqn:Z2LinComb}), the only tangent vector having non-zero $(i_{n,1},j_n)$-component is $v_{p,n,1}$; whence, the coefficient $d_{n,1}$ $=$ $0$.  Since $v_{p,n,1}$ is the only tangent vector in $\EuScript{Z}'$ associated to $\alpha$ $=$ $n$, we have shown that all the coefficients $d_{n,\delta}$ are $0$; this completes the base case of the induction.  

For the induction step, we suppose that for some $1 \leq \alpha < n$, the coefficients $d_{\beta,\delta'}$ $=$ $0$ for all $\beta$ $>$ $\alpha$ and $1 \leq \delta' \leq |\Delta'_{\beta}|$, and let $1 \leq \delta \leq |\Delta'_{\alpha}|$.  We claim that the $(i_{\alpha,\delta},j_{\alpha})$-component of $v_{p,\beta,\delta'}$ is $0$ for all $\beta$ $\leq$ $\alpha$ and $(\beta,\delta')$ $\neq$ $(\alpha,\delta)$.  This component is the coefficient of $t_{i_{\alpha,\delta}}$ in $\left( \frac{\partial g_{j_{\alpha}}}{\partial x_{\beta}} \cdot m_{\beta,\delta'} \mod I_p\right)$. In case $b_{j_{\alpha}}$ $\notin$ $\leadmon$, we have $g_{j_{\alpha}}$ $=$ $b_{j_{\alpha}}$ $=$ $x_{\alpha} x_n^{e'_{\alpha}}$; otherwise, $b_{j_{\alpha}}$ $\in$ $\leadmon$ and $g_{j_{\alpha}}$ $=$ $b_{j_{\alpha}}- N_{j_{\alpha}}$, and by hypothesis $b_{j_{\alpha}}$ $>$ any monomial $m$ appearing non-trivially in $N_{j_{\alpha}}$.  From this it follows that for all $\beta$ $<$ $\alpha$, $\frac{\partial g_{j_{\alpha}}}{\partial x_{\beta}}$ $=$ $0$.  This shows that the $(i_{\alpha,\delta},j_{\alpha})$-coefficient of $v_{p,\beta,\delta'}$ is $0$ for all $\beta$ $<$ $\alpha$.

We now consider the case $\beta$ $=$ $\alpha$ and $\delta'$ $\neq$ $\delta$.  We must compute the coefficient of $t_{i_{\alpha,\delta}}$ in $\left( \frac{\partial g_{j_{\alpha}}}{\partial x_{\alpha}} \cdot m_{\alpha,\delta'} \mod I_p \right)$.  In light of $(\ref{txt:keyClaim})$, we see that $\frac{\partial g_{j_{\alpha}}}{\partial x_{\alpha}}$ $=$ $\frac{\partial b_{j_{\alpha}}}{\partial x_{\alpha}}$ $=$ $x_n^{e'_{\alpha}}$, so $\left(\frac{\partial g_{j_{\alpha}}}{\partial x_{\alpha}} \cdot m_{\alpha,\delta'} \mod I_p \right)$ = $(t_{i_{\alpha,\delta'}} \mod I_p)$ $=$ $t_{i_{\alpha,\delta'}}$ $\neq$ $t_{i_{\alpha,\delta}}$ for $\delta$ $\neq $ $\delta'$.  This shows that the $(i_{\alpha,\delta},j_{\alpha})$-component of $v_{p,\beta,\delta'}$ is $0$ for all $\beta$ $\leq$ $\alpha$ and $(\beta,\delta')$ $\neq$ $(\alpha,\delta)$. It now follows from (\ref{txt:nonSpecialIndices}) and the induction hypothesis that   
the coefficients $d_{\alpha,\delta}$ $=$ $0$ for all $1 \leq \delta \leq |\Delta'_{\alpha}|$, so the induction step is complete. We conclude that all the coefficients $d_{\alpha,\delta}$ in (\ref{eqn:Z2LinComb}) are $0$.  Lemma \ref{lem:YtanVecsLinInd} now yields that every scalar $d_{ij}$ $=$ $0$ as well, and the proposition is proved.
\qed

\begin{cor} \label{cor:GenZPropCor1A}
    Under the hypotheses of Proposition {\rm \ref{prop:GeneralZProp}}, one has that the $(i_{\alpha,\delta},j_{\alpha})$-component of $v_{p,\beta,\delta'}$ $=$ $0$ for all $\beta$ $\leq$ $\alpha$ and $(\beta,\delta')$ $\neq$ $(\alpha,\delta)$.  
\end{cor}

\proof
    This was shown in the course of the proof of the Proposition.
\qed
   
\begin{cor}\label{cor:GenZPropCor1}
    Let $\mathcal{O}$ $\neq$ $\{ 1 \}$ be a lex-segment complement order ideal.  Then the set $\EuScript{S}_0 \cup \EuScript{Z'}_0$ $\subseteq$ $\tansp_{[I_0]}$ is $\gf$-linearly independent.     
(Recall that $I_0$ is the \SLI\ ideal generated by the monomials $m$ $\notin$ $\mathcal{O}$.) 
\end{cor}

\proof
    Since $b_{j_{\alpha}}$ is the lex-leading monomial of $g_{j_{\alpha}}$ $=$ $b_{j_{\alpha}}$ for $1$ $\leq$ $\alpha$ $\leq$ $n$, Proposition {\rm \ref{prop:GeneralZProp}} yields the result.
\qed

\begin{cor} \label{cor:GenZPropCor2}
    Let $\mathcal{O}$ $\neq$ $\{ 1 \}$ be a lex-segment complement order ideal, and suppose that the sets $\leadmon$ and $\trailmon$ have been chosen such that for all $m \in \leadmon$ and $m' \in \trailmon$, one has that $m>m'$.  Then the  set $\EuScript{S}_p \cup \EuScript{Z'}_p$ $\subseteq$ $\tansp_{[I_p]}$ is $\gf$-linearly independent for all points $p$ as in {\rm (\ref{eqn:pRestriction})}; that is, for all \SLI\ ideals $I$ $=$ $I_p$.
\end{cor}

\proof
    The hypotheses clearly imply that the hypotheses of Proposition \ref{prop:GeneralZProp} hold for all $p$; whence, the result.
\qed

\begin{cor} \label{cor:GenZPropCor3}
    Let $\mathcal{O}$ $\neq$ $\{ 1 \}$ be a lex-segment complement order ideal, let $U$ $=$ $\Spec(\gf[\mathcal{C}, \mathcal{Z}'])$, and $\fn{\EuScript{F}}{U}{\Hilb^n_{\mathbb{A}^n_{\gf}}}$ the map constructed in Section {\rm \ref{ssec:Cconstr}}.  Then the \SLI\ locus $X_{\specl}$ is contained in a component $Y$ of $\Hilb^{\mu}_{\mathbb{A}^n_{\gf}}$ of dimension
\[
    \dim(Y) \geq  \dim(\EuScript{F}(U)) = \dim(U) = |\specl| + |\mathcal{Z}'|.
\]
\end{cor} 

\proof
    By Proposition \ref{prop:F(U)DimLowerBd} and Corollary \ref{cor:GenZPropCor1}, we have that
\[
    |\specl| + |\mathcal{Z}'| = |\EuScript{S}_0| + |\EuScript{Z}'_0| \leq \dim(\EuScript{F}(U)) \leq \dim(U) = |\specl| + |\mathcal{Z}'|.
\]
Choosing  $Y$ to be a component of $\Hilb^{\mu}_{\mathbb{A}^n_{\gf}}$ containing $\EuScript{F}(U)$, we have that $Y$ contains $X_{\specl}$ and $\dim{Y}$ $\geq$ $\dim(\EuScript{F}(U))$ $=$ $|\specl| + |\mathcal{Z}'|$, as desired. 
\qed

\subsection{Additional independent tangent directions at $[I_p]$} \label{ssec:ZExtensionRem}
Proposition \ref{prop:GeneralZProp} and its corollaries give conditions on the \SLI\ ideal $I_p$ for which the set $\EuScript{S}_p$ $\cup$ $\EuScript{Z}'_p$ $\subseteq$ $\tansp_{[I_p]}$ is $\gf$-linearly independent.  For many of our examples this suffices, because in these cases $\EuScript{S}_p$ $\cup$ $\EuScript{Z}'_p$ turns out to be a $\gf$-basis of $\tansp_{[I_p]}$.  However, in other cases $\tansp_{[I_p]}$ has a basis that is a proper superset of $\EuScript{S}_p$ $\cup$ $\EuScript{Z}'_p$; here is one way this can happen.

With $\mathcal{O}$ $\neq$ $\{ 1 \}$ a lex-segment complement, suppose that we have chosen the sets $\Delta_{\alpha}$ $\subseteq$ $\gf[x_1, \dots, \widehat{x_{\alpha}},\dots, x_n]$ to be supersets of the corresponding sets $\Delta'_{\alpha}$.  Furthermore, suppose that at least one of the monomials $b_{j_{\alpha}}$ is non-leading and that there is a monomial $m_{\alpha,\hat{\delta}}$ $\in$ $\Delta_{\alpha}$ such that $t_{i_{\alpha}} \cdot m_{\alpha,\hat{\delta}}$ $=$ $b_{\hat{j}}$ $\in$ $\leadmon$, with
\[
     g_{\hat{j}} = \left( b_{\hat{j}} - \sum_{C_{i\hat{j}} \in \specl} {c_{i\hat{j}}\, t_{i}}\right)\ \neq\ b_{\hat{j}}.
\]
  Then the non-zero $(i,j_{\alpha})$-components of the tangent vector $v_{\alpha,\hat{\delta}}$ are the (non-zero) coefficients $c_{i\hat{j}}$ of the linear combination
\[
     t_{i_{\alpha}} \cdot m_{\alpha,\hat{\delta}} = b_{\hat{j}} \equiv \sum_{C_{i\hat{j}} \in \specl} {c_{i\hat{j}}\, t_{i}}\ \text{ mod } I,
\]   
and these components are non-\SLI\ because $b_{j_{\alpha}}$ $\notin$ $\leadmon$.  It is therefore possible that $v_{p,\alpha,\hat{\delta}}$ $\in$ $\tansp_{[I]}$ is independent of the vectors in $\EuScript{Z}'$ (and certain that it is independent of the vectors in $\EuScript{S}$).   Accordingly, we define
\[
    \begin{array}{rcl}
        \Delta''_{\alpha} & = & \Delta'_{\alpha} \cup \left\{ m_{\alpha,\hat{\delta}} \in \gf[x_1, \dots, \widehat{x_{\alpha}}, \dots, x_n] \mid 
        \begin{array}{r}
               b_{j_{\alpha}} \notin \leadmon \text{ and }\\      
                t_{i_{\alpha}} \cdot m_{\alpha,\hat{\delta}} \in \leadmon
        \end{array} \right\},\vspace{.05in}\\
 \mathcal{Z}'' & = & \{ Z_{\alpha,\delta} \mid 1 \leq \alpha \leq n,\ 1 \leq \delta \leq |\Delta''_{\alpha}| \},\vspace{.05in}\\
        \EuScript{Z}_p'' & = & \{v_{p,\alpha,\delta} \mid 1 \leq \alpha \leq n,\ 1 \leq \delta \leq |\Delta''_{\alpha}| \} \subseteq \tansp_{[I_{p}]} .
    \end{array} 
\]

From Proposition \ref{prop:F(U)DimLowerBd}, we obtain the following analogue of Corollary \ref{cor:GenZPropCor3}:

\begin{cor} \label{cor:Z''Cor}
    Let $\mathcal{O}$ $\neq$ $\{ 1 \}$ be a lex-segment complement order ideal, let $U$ $=$ $\Spec(\gf[\mathcal{C}, \mathcal{Z}''])$, and $\fn{\EuScript{F}}{U}{\Hilb^n_{\mathbb{A}^n_{\gf}}}$ the map constructed in Section {\rm \ref{ssec:Cconstr}}.  Suppose that there is a \SLI\ ideal $I_p$ such that the set $\EuScript{S}_p$ $\cup$ $\EuScript{Z}''_p$ $\subseteq$ $\tansp_{[I_p]}$ is $\gf$-linearly independent. Then the \SLI\ locus $X_{\specl}$ is contained in a component $Y$ of $\Hilb^{\mu}_{\mathbb{A}^n_{\gf}}$ of dimension
\[
    \dim(Y) \geq  \dim(\EuScript{F}(U)) = \dim(U) = |\specl| + |\mathcal{Z}''|.
\] 
\qed
\end{cor}

In several of the following examples (Sections \ref{ssec:(1,5,3,4)}, \ref{ssec:(1,5,3,4,5,6)}, \ref{ssec:(1,6,6,10)}, \ref{ssec:(1,6,21,10,15)}, and \ref{(1,6,10,10,5)CaseA}), we have that the set $\EuScript{S}_p \cup \EuScript{Z}'_p$ is a $\gf$-basis of $\tansp_{[I_p]}$, so expanding to $\EuScript{Z}''$ gives us nothing new.  On the other hand, the examples presented in Sections \ref{(1,6,10,10,5)CaseB} and \ref{(1,6,10,10,5)CaseC} are such that $\EuScript{S}_p \cup \EuScript{Z}'_p$ is not a basis of $\tansp_{[I_p]}$, but $\EuScript{S}_p \cup \EuScript{Z}''_p$ is.


\section{Examples of generic \SLI\ ideals} \label{sec:examples}

	We present several examples of generic \SLI\ ideals associated to lex-segment complement order ideals.  Unfortunately, the examples are too large to permit the computations to be carried out by hand. We summarize each example and provide the details in a \textit{Mathematica} \cite{Mathematica} notebook that is available for download from the arXiv, where it is posted as an ancillary file to this paper.  The notebooks are also available at
\[
    \text{\emph{http://www.skidmore.edu/\~{}mhuibreg/Notebooks for paper/index.html}}.
\] 
The notation used in the notebooks adheres closely to that used in the paper.  The notebooks all make use of a library of \textit{Mathematica} functions coded and documented in a separate notebook \emph{utility functions.nb} that is available for download at the same locations.  

	Each example of a generic ideal is a \SLI\ ideal $I$ such that $[I]$ is a smooth point on an elementary  component of the form $\overline{\EuScript{F}(U)}$ as in Proposition \ref{prop:genericCond}. The ideal $I$ is generated by a \emph{Mathematica} function \textbf{makeShortListOfIdealGenerators}, which, given the sets $\leadmon$ and $\trailmon$, generates the polynomials $g_{j_{\iota}}$ $\in$ $G$ (equation \ref{eqn:fDef}) by assigning values drawn at random (with equal probabilities) from the set $\{-1,0,1 \}$ to the \SLI\ coefficients $C_{ij}$. (A second version of this function is provided that assigns the coefficients from the set $\{0,1 \}$ with the probability of assigning the value $1$ supplied as an additional input; this version was used to generate the example discussed in Remark \ref{rem:GenButNotEff}, and nowhere else.)  
Since $\overline{\EuScript{F}(U)}$ will be smooth in a neighborhood of a smooth point, the existence of a generic \SLI\ ideal $I$ implies that there is a non-empty Zariski-open subset $X'$ $\subseteq$ $X_{\specl}$ such that $[I']$ $\in$ $X'$ $\Rightarrow$ $I'$ is generic.  In each example notebook, one can either choose to verify that a previously-generated example (stored in the notebook) is generic, or can generate other \SLI\ ideals to test.  In view of the foregoing, such examples will typically (but need not always) be generic as well.  Note that in all cases the order ideal $\mathcal{O}$ is the unique lex-segment complement having the given Hilbert function.


\subsection{Hilbert function $(1,5,3,4,0)$, shape $(5,2,2,3)$} \label{ssec:(1,5,3,4)} 

	We exhibit a generic \SLI\ ideal $I$ $=$ $I_p$ of this Hilbert function in the notebook \emph{case $(1,5,3,4,0)$.nb}.  For this initial example, we provide more details here to serve as an introduction; recall that this example was summarized in Section \ref{ssec:AGen}. 

	In this example, we have 
\[
  \begin{array}{l}
    n  =  5, \vspace{.05in}\\
     \mathcal{O}  =  \left\{1,x_1,x_2,x_3,x_4,x_5,x_4^2,x_4x_5,x_5^2,x_4^3,x_4^2x_5,x_4x_5^2,x_5^3\right\};\vspace{.05in}\\
\leadmon  =  \left\{x_1^2,x_1 x_2,x_1 x_3,x_1 x_4,x_1 x_5,x_2^2,x_2 x_3,x_2 x_4,x_2 x_5,x_3^2,x_3 x_4,x_3 x_5\right\};\vspace{.05in}\\
\trailmon  =  \left\{x_4^3,x_4^2 x_5,x_4 x_5^2,x_5^3\right\};\vspace{.05in}\\
\lambda = |\leadmon| = 12,\ \tau = |\trailmon| = 4,\ 
\mu = |\mathcal{O}| = 13.
  \end{array}
\]

As in Section \ref{sec:SLI}, we constructed the following set of polynomials $G$ $=$ \(\left\{\left.g_j\right|1\leq j\leq \lambda \right\}\), and extended them to the $\mathcal{O}$-border basis of a \SLI\  ideal \(I\) that can then be shown to be efficient (in fact, $\vartheta$-efficient). 
{\Small
\[ 
\left\{
  \begin{array}{ll}
     g_1=x_1^2-x_4^3-x_4^2 x_5+x_4 x_5^2+x_5^3, & g_2=x_1 x_2-x_4^3+x_4^2 x_5+x_4 x_5^2+x_5^3, \vspace{.05in}\\ 
g_3=x_1 x_3-x_4^2 x_5+x_5^3, &
    g_4=x_1 x_4+x_4^3+x_4^2 x_5-x_5^3,\vspace{.05in}\\
g_5=x_1 x_5 + x_4^2 x_5-x_5^3, & g_6=x_2^2+x_4^3+x_4^2 x_5+x_4 x_5^2-x_5^3,\vspace{.05in}\\
    g_7=x_2 x_3-x_4^3-x_4^2 x_5, & g_8=x_2 x_4+x_4^3-x_4^2 x_5,\vspace{.05in}\\
 g_9=x_2 x_5-x_4^2 x_5, & g_{10}=x_3^2+x_4^3+x_4^2 x_5+x_4
x_5^2-x_5^3,\vspace{.05in}\\
 g_{11}=x_3 x_4+x_4^3-x_4^2 x_5, & g_{12}=x_3 x_5-x_4^2 x_5-x_4 x_5^2 
  \end{array} 
\right\}
\]
}

We list the sets of monomials $\Delta'_{\alpha}$ as in (\ref{eqn:DeltaDelta'RestrDef}):
\[
    \begin{array}{cccc}
        x_{\alpha} & b_{j_{\alpha}} & t_{i_{\alpha}} & \Delta'_{\alpha}\vspace{.05in}\\
         x_1  & x_1 x_5 \in \leadmon & x_5 & \{ 1, x_4, x_5 \}\vspace{.05in}\\
         x_2  & x_2 x_5 \in \leadmon & x_5 & \{ 1, x_4, x_5 \}\vspace{.05in}\\
         x_3  & x_3 x_5 \in \leadmon & x_5 & \{ 1, x_4, x_5 \}\vspace{.05in}\\
         x_4  & x_4 x_5^3 \notin \leadmon & x_5^3 & \{ 1 \}\vspace{.05in}\\
         x_5  & x_5^4 \notin \leadmon & x_5^3 & \{ 1 \}
    \end{array}
\]

By Corollary \ref{cor:GenZPropCor3} the image of the associated map $\fn{\EuScript{F}}{U}{\Hilb^{13}_{\mathbb{A}^5_{\gf}}}$ satisfies
\[
    \dim(\EuScript{F}(U)) \geq |\specl| + |\mathcal{Z}'| = (\lambda \cdot \tau) + (3+3+3+1+1)  = 12 \cdot 4 + 11 = 59.
\]
On the other hand, when we compute the dimension of the tangent space $\tansp_{[I_p]}$ using the tangent space relations associated to a basis of linear syzygies, as in (\ref{eqn:tanSpDimComp}), we obtain that $\dim_{\gf}(\tansp_{[I_p]})$ $=$ $59$; consequently, by Proposition \ref{prop:genericCond}, $\overline{\EuScript{F}(U)}$ is an elementary component of $\Hilb^{13}_{\mathbb{A}^5_{\gf}}$ of dimension  59, on which $[I_p]$ is a smooth point.  Note that the dimension of the principal component is $5 \cdot 13$ $=$ $65$.

\begin{rem} \label{rem:GenButNotEff}
    The notebook also includes an example of a generic \SLI\ ideal with Hilbert function $(1,5,3,4,0)$ that is not efficient.  Therefore, efficiency is not necessary for genericity.
\end{rem}


\subsection{Hilbert function $(1,5,3,4,5,6,0)$, shape $(5,2,2,5)$} \label{ssec:(1,5,3,4,5,6)} 

	The details of this example are presented in the notebook \emph{case $(1,5,3,4,5,$ $6,0)$.nb}.  The set of leading monomials is the same as in the previous example, and the set of trailing monomials consists of the six monomials of degree $5$ in $x_4, x_5$.  Sufficiently general \SLI\ ideals $I$ of this shape are $\vartheta$-efficient.

	The sets $\Delta'_{\alpha}$ are computed as in the previous example:  
\[
    \begin{array}{c}
    \Delta'_1 = \Delta'_2 = \Delta'_3 = \{ 1,\, x_5,\, x_6,\, x_5^2,\, x_5 x_6,\, x_6^2,\,x_5^3,\,x_5^2 x_6,\, x_5 x_6^2,\,x_6^3\}, \text{ and }\vspace{.05in}\\
         \Delta'_4 = \Delta'_5  = \{1\}
    \end{array}.
\]
By Corollary \ref{cor:GenZPropCor3} the image of the associated map $\fn{\EuScript{F}}{U}{\Hilb^{24}_{\gf}}$ satisfies
\[
    \dim(\EuScript{F}(U)) = L \geq |\specl| + |\mathcal{Z}'| = \lambda \cdot \tau + 3\cdot 10 + 2 = 12 \cdot 6 + 32 = 104.
\]
Moreover, for sufficiently general choices of the generators $G$, we have that $\dim_{\gf}(\tansp_{[I]})$ $=$ $104$; consequently, by Proposition \ref{prop:genericCond}, $\overline{\EuScript{F}(U)}$ is an elementary component of $\Hilb^{24}_{\mathbb{A}^5_{\gf}}$ of dimension 104, on which $[I]$ is a smooth point.  Note that the dimension of the principal component is $5 \cdot 24$ $=$ $120$ in this case. 


\subsection{Hilbert function $(1,6,6,10,0)$, shape $(6,3,2,3)$} \label{ssec:(1,6,6,10)}

The details of this example are presented in the notebook \emph{case $(1,6,6,10,0)$.nb}.  
The sets $\leadmon$, $\trailmon$, and $\Delta'_{\alpha}$ are as follows:
\[
    \begin{array}{c}
        \leadmon\  =\  \left\{x_1^2,\, x_1 x_2,\, \dots,\, x_3 x_6\right\},\ \   \trailmon\  =\  \left\{ x_4^3,\, x_4^2 x_5,\, \dots,\, x_6^3  \right\},\vspace{.05in}\\        
      \Delta'_1 = \Delta'_2 = \Delta'_3 = \{ 1,\, x_4,\, x_5,\, x_6 \},\ \ \Delta'_4 = \Delta'_5 = \Delta'_6 = \{ 1 \}.
    \end{array}
\]
By Corollary \ref{cor:GenZPropCor3} the image of the associated map $\fn{\EuScript{F}}{U}{\Hilb^{23}_{\gf}}$ satisfies
\[
    \dim(\EuScript{F}(U)) = L \geq |\specl| + |\mathcal{Z}'| = \lambda \cdot \tau + 3\cdot 4 + 3 = 15\cdot 10  + 15 = 165.
\]
Moreover, for sufficiently general choices of the set $G$, we have that $\dim_{\gf}(\tansp_{[I]})$ $=$ $165$; consequently, by Proposition \ref{prop:genericCond}, $\overline{\EuScript{F}(U)}$ is an elementary component of $\Hilb^{23}_{\mathbb{A}^5_{\gf}}$ of dimension 165, on which $[I]$ is a smooth point.   Note that the dimension of the principal component is $6 \cdot 23$ $=$ $138$.
  
\begin{rem} \label{rem:genButNotEff} 
     In the notebook detailing this example, we observe that sufficiently general \SLI\ ideals are efficient, but no \SLI\ ideal in this case can be $\vartheta$-efficient.  Consequently, $\vartheta$-efficiency is not a necessary condition for genericity,  as we noted in Remark \ref{rem:ThEffAndEff}. 
\end{rem}


\subsection{Hilbert function $(1,6,21,10,15,0)$, shape $(6,3,3,4)$} \label{ssec:(1,6,21,10,15)}   

The details of this example are presented in the notebook \emph{case $(1,6,21,10,15,$ $0)$.nb}.  The set of leading monomials is equal to the first 46 monomials of degree $3$ when listed in decreasing lex order, and the set of trailing monomials consists of the 15 monomials of degree $4$ in $x_4, x_5, x_6$.  One finds that sufficiently general \SLI\ ideals $I$ of this shape are $\vartheta$-efficient.  The sets $\Delta'_{\alpha}$ are computed as in (\ref{eqn:DeltaDelta'RestrDef}): 
\[
    \Delta'_1 = \Delta'_2 = \Delta'_3 = \{ 1, x_4, x_5, x_6\}, \text{ and } \Delta'_4 = \Delta'_5 = \Delta'_6 = \{1\}.
\]
From Corollary \ref{cor:GenZPropCor3} we obtain
\[
   \dim(\EuScript{F}(U)) \geq |\specl| + |\mathcal{Z}'| = 46 \cdot 15 + 4\cdot 3 +3  =  705;
\]
moreover, Corollary \ref{cor:GenZPropCor2} implies that the tangent space at $[I]$ for \emph{every} \SLI\ ideal $I$ has dimension $\geq 705$.  The number of variables $a_{ij}$ involved in the tangent space relations (Section \ref{sec:tanSpace}) is
\[
    (\text{\# boundary monoms}) \cdot (\text{\# basis monomials}) =  142\cdot 53 = 7526. 
\] 
From this it follows that  the rank $\rho$ of the tangent space relations at 
$[I]$ satisfies
\[
    7526 - \rho \geq 705\ \Rightarrow\ 6821 \geq \rho.
\] 
It follows that $I$ will be generic provided that the rank of the tangent space relations at $[I]$ is equal to its maximum possible value of $6821$.  
In the notebook associated to this example, we computed this rank modulo a large prime (32713) to conserve memory, and obtained the value $6821$.  Since the tangent space relations in our examples have integer coefficients, and the rank of an integer matrix cannot increase when one computes it modulo a prime, this computation demonstrates that the characteristic $0$ rank must be $6821$.     
Therefore, $\dim_{\gf}(\tansp_{[I]})$ $=$ $7526 - 6821$ $=$ $705$, so, once again, Proposition \ref{prop:genericCond} implies that $[I]$ is a smooth point on an elementary component of $\Hilb^{53}_{\mathbb{A}^6_{\gf}}$ of dimension $705$.  The dimension of the principal component is $6 \cdot 53$ $=$ $318$.   Note that the computations in this notebook require a lot of memory (two gigabytes was insufficient) and a run time possibly measured in hours, depending on the speed of one's machine.


\subsection{Hilbert function $(1,6,10,10,5,0)$} \label{ssec:(1,6,10,10,5)}

     We present three different examples of elementary components having the indicated Hilbert function.  The order ideal for all three examples is the lex-segment complement
\[
\mathcal{O} = 
    \left\{ \begin{array}{c}
   1\vspace{.05in}\\
x_1\ \  x_2\ \ x_3\ \ x_4\ \ x_5\ \ x_6 \vspace{.05in}\\
x_3^2\ \ x_3 x_4\ \ x_3 x_5\ \ x_3 x_6\ \ x_4^2\ \ x_4 x_5\ \ x_4 x_6\ \ x_5^2\ \ x_5 x_6\ \ x_6^2 \vspace{.05in}\\
x_4^3\ \ x_4^2 x_5\ \ x_4^2 x_6\ \ x_4 x_5^2\ \ x_4 x_5 x_6\ \ x_4 x_6^2\ \ x_5^3\ \ x_5^2 x_6\ \ x_5 x_6^2\ \ x_6^3 \vspace{.05in}\\x_5^4\ \  x_5^3 x_6\ \   x_5^2 x_6^2\ \ x_5 x_6^3\ \ x_6^4 
    \end{array} \right\}
\]  
of cardinality $32$.  Note that the dimension of the principal component of $\Hilb^{32}_{\mathbb{A}^6_{\gf}}$ is $32 \cdot 6$ $=$ $192$.

\subsubsection{First case} \label{(1,6,10,10,5)CaseA}

    The details of this example are presented in the notebook \emph{case $(1,6,10,10,5,0)$ first.nb}.
In this case the sets of leading and trailing monomials are
\[
    \begin{array}{rcl}
        \leadmon & = & \left\{ \begin{array}{c} 
                                    x_1^2,\,x_1 x_2,\,x_1 x_3,\,x_1 x_4,\,x_1 x_5,\,x_1 x_6,\,x_2^2,\vspace{.05in}\\ 
x_2 x_3,\,x_2 x_4,\,x_2 x_5,\,x_2 x_6, x_3^3,\,x_3^2 x_4,\,x_3^2 x_5,\, x_3^2 x_6,
\vspace{.05in}\\
x_3 x_4^2,\,x_3 x_4 x_5,\,x_3 x_4 x_6,\,x_3 x_5^2,\,x_3 x_5 x_6,\,x_3 x_6^2
                              \end{array} \right\}, \vspace{.05in}\\
\trailmon & = & \left\{ \begin{array}{c}x_4^3,\, x_4^2 x_5, \,x_4^2 x_6, \,x_4 x_5^2, \,x_4 x_5 x_6,\vspace{.05in}\\
   x_4 x_6^2,\,x_5^4,\,x_5^3 x_6,\,x_5^2 x_6^2,\,x_5 x_6^3,\,x_6^4
                        \end{array} \right\}.
    \end{array}
\]
One finds that sufficiently general \SLI\ ideals $I$ constructed using these sets are efficient.  

The sets $\Delta'_{\alpha}$ are as follows, showing that $|\mathcal{Z}'|$ $=$ $24$:
\[ 
    \begin{array}{cccc}
        x_{\alpha} & b_{j_{\alpha}} & t_{i_{\alpha}} & \Delta'_{\alpha}\vspace{.05in}\\
         x_1  & x_1 x_6 \in \leadmon & x_6 & \{ 1,\, x_3,\, x_4,\, x_5,\, x_6,\, x_5^2,\, x_5 x_6,\, x_6^2 \}\vspace{.05in}\\
         x_2  & x_2 x_6 \in \leadmon & x_6 & \{ 1,\, x_3,\, x_4,\, x_5,\, x_6,\, x_5^2,\, x_5 x_6,\, x_6^2 \}\vspace{.05in}\\
         x_3  & x_3 x_6^2 \in \leadmon & x_6^2 & \{ 1,\, x_5,\, x_6 \}\vspace{.05in}\\
         x_4  & x_4 x_6^3 \notin \leadmon & x_6^3 & \{ 1,\, x_5,\, x_6 \}\vspace{.05in}\\
         x_5  & x_5 x_6^4 \notin \leadmon & x_6^4 & \{ 1 \}\vspace{.05in}\\
         x_6  & x_6^5 \notin \leadmon & x_6^4 & \{ 1 \}
    \end{array}.
\]
From Corollary \ref{cor:GenZPropCor3} we obtain that
\[
    \dim(\EuScript{F}(U)) \geq |\specl| + |\mathcal{Z}'| = \lambda \cdot \tau + 24 = 21 \cdot 11 + 24 = 255.
\]
On the other hand, we find by direct computation in the notebook that $\dim_{\gf}(\tansp_{[I]})$ $=$ $255$, so Proposition \ref{prop:genericCond} implies that $[I]$ is a smooth point on an elementary component of $\Hilb^{32}_{\mathbb{A}^6_{\gf}}$ of dimension $255$.


\subsubsection{Second case} \label{(1,6,10,10,5)CaseB}

     The details of this example are presented in the notebook \emph{case $(1,6,10,10,5,0)$ second.nb}. The sets of leading and trailing monomials are
\[
    \begin{array}{rcl}
        \leadmon & = & \left\{ \begin{array}{c} 
                                    x_1^2,\, x_1 x_2,\, x_1 x_3,\, x_1 x_4,\, x_1 x_5,\, x_1 x_6,\, x_2^2,\, x_2 x_3,\vspace{.05in}\\
 x_2 x_4,\, x_2 x_5,\,x_2 x_6,\, x_4^4,\, x_4^3 x_5,\, x_4^3 x_6,\,x_4^2 x_5^2,\vspace{.05in}\\x_4^2 x_5 x_6,\,x_4^2 x_6^2,\,x_4 x_5^3,\,x_4 x_5^2 x_6,\,x_4 x_5 x_6^2,\,x_4 x_6^3
                              \end{array} \right\}, \vspace{.05in}\\
\trailmon & = & \left\{x_3^2,\,x_3 x_4,\,x_3 x_5,\,x_3 x_6,\,x_5^4,\,x_5^3 x_6,\,x_5^2 x_6^2,\,x_5 x_6^3,\,x_6^4\right\}. 
    \end{array}
\]
As usual, a \SLI\ ideal $I$ $=$ $I_p$ is generated, and its tangent space dimension is computed to be $\dim_{\gf}(\tansp_{[I]})$ $=$ $222$; the ideal $I$ is also found to be efficient.  

In this case, there are trailing monomials that are lex-larger than some leading monomials; for example, $x_3 x_6$ $>$ $x_4 x_6^3$ $=$ $b_{j_4}$, so it is likely that \SLI\ ideals $[I]$ will fail to satisfy the hypothesis of Proposition \ref{prop:GeneralZProp}. On the other hand, there is a non-leading boundary monomial $b_{j_3}$ $=$ $x_3 x_6^2$ such that $t_{i_3}\cdot x_4 x_6$ $=$ $x_4 x_6^3$ $=$ $b_{j_4}$ $\in$ $\leadmon$, so the situation described in Section \ref{ssec:ZExtensionRem} arises; that is, it is possible that the larger set of tangent vectors $\EuScript{S}_p \cup \EuScript{Z}''_p$  is $\gf$-linearly independent at $[I_p]$.  Indeed, we verify this by direct computation in the notebook. The associated sets $\Delta''_{\alpha}$ are as follows: 
\[ 
    \begin{array}{cccc}
        x_{\alpha} & b_{j_{\alpha}} & t_{i_{\alpha}} & \Delta''_{\alpha}\vspace{.05in}\\
         x_1  & x_1 x_6 \in \leadmon & x_6 & \{ 1,\, x_4,\, x_5,\, x_6,\, x_4^2,\, x_4 x_5,\, x_4 x_6,\, x_5^2,\, x_5 x_6,\, x_6^2 \}\vspace{.05in}\\
         x_2  & x_2 x_6 \in \leadmon & x_6 & \{ 1,\, x_4,\, x_5,\, x_6,\, x_4^2,\, x_4 x_5,\, x_4 x_6,\, x_5^2,\, x_5 x_6,\, x_6^2 \}\vspace{.05in}\\
         x_3  & x_3 x_6^2 \notin \leadmon & x_6^2 & \{ 1,\, x_4,\, x_5,\, x_6,\, x_4^2,\, x_4 x_5,\, x_4 x_6,\, x_5^2,\, x_5 x_6,\, x_6^2 \}\vspace{.05in}\\
         x_4  & x_4 x_6^3 \in \leadmon & x_6^3 & \{ 1 \}\vspace{.05in}\\
         x_5  & x_5 x_6^4 \notin \leadmon & x_6^4 & \{ 1 \}\vspace{.05in}\\
         x_6  & x_6^5 \notin \leadmon & x_6^4 & \{ 1 \}
    \end{array}.
\]
By Corollary \ref{cor:Z''Cor} the image of the map $\fn{\EuScript{F}}{U}{\Hilb^{32}_{\mathbb{A}^6_{\gf}}}$ satisfies
\[
    \dim(\EuScript{F}(U)) \geq |\specl| + |\mathcal{Z}''| = (\lambda \cdot \tau) + (10 \cdot 3 + 3)  = 21 \cdot 9 + 33 = 222.
\]
It now follows from Proposition \ref{prop:genericCond} that $[I]$ is a smooth point on an elementary component of dimension $222$, so that $I$ is a generic ideal.


\subsubsection{Third case} \label{(1,6,10,10,5)CaseC}
 The details of this example are presented in the notebook \emph{case $(1,6,10,10,5,0)$ third.nb}. The sets of leading and trailing monomials are
\[
    \begin{array}{rcl}
        \leadmon & = & \left\{ \begin{array}{c} 
                                  x_1^2,\, x_1 x_2,\, x_1 x_3,\, x_1 x_4,\, x_1 x_5,\, x_1 x_6,\, x_2^2,\, x_2 x_3,\vspace{.05in}\\ x_2 x_4,\, x_2 x_5,\, x_2 x_6,
x_5^5,\, x_5^4 x_6,\, x_5^3 x_6^2,\, x_5^2 x_6^3,\, x_5 x_6^4,\, x_6^5
                              \end{array} \right\} \vspace{.05in}\\
\trailmon & = & \left\{ \begin{array}{c} x_3^2,\, x_3 x_4,\, x_3 x_5,\, x_3 x_6,\, x_4^3,\, x_4^2 x_5,\vspace{.05in}\\ x_4^2 x_6,\, x_4 x_5^2,\, x_4 x_5 x_6,\, 
x_4 x_6^2
                        \end{array} \right\}
    \end{array}.
\]

We randomly generate a \SLI\ ideal $I$ $=$ $I_p$ using these sets and find that $\dim_{\gf}(\tansp_{[I]})$ $=$ $211$; we also find that $I$ is efficient, but not $\vartheta$-efficient, thereby providing another example as promised in Remark
\ref{rem:ThEffAndEff}. 

As in the preceding example, we verify by direct computation that the larger set of vectors $\EuScript{S}_p \cup \EuScript{Z}''_p$ is $\gf$-linearly independent. The associated sets $\Delta''_{\alpha}$ are as follows:
{\Small 
\[ 
    \begin{array}{cccc}
        x_{\alpha} & b_{j_{\alpha}} & t_{i_{\alpha}} & \Delta''_{\alpha}\vspace{.05in}\\
         x_1  & x_1 x_6 \in \leadmon & x_6 & \{ 1,\, x_4,\, x_5,\, x_6,\, x_5^2,\, x_5 x_6,\,  x_6^2,\, x_5^3,\, x_5^2 x_6,\,x_5 x_6^2, x_6^3 \}\vspace{.05in}\\
         x_2  & x_2 x_6 \in \leadmon & x_6 & \{ 1,\, x_4,\, x_5,\, x_6,\, x_5^2,\, x_5 x_6,\,  x_6^2,\, x_5^3,\, x_5^2 x_6,\,x_5 x_6^2, x_6^3 \}\vspace{.05in}\\
         x_3  & x_3 x_6^2 \notin \leadmon & x_6^2 & \{1,\, x_4,\, x_5,\, x_6,\, x_5^2,\, x_5 x_6,\,  x_6^2,\, x_5^3,\, x_5^2 x_6,\,x_5 x_6^2, x_6^3 \}\vspace{.05in}\\
         x_4  & x_4 x_6^3 \notin \leadmon & x_6^3 & \{ 1,\, x_5,\, x_6,\, x_5^2,\, x_5 x_6,\, x_6^2 \}\vspace{.05in}\\
         x_5  & x_5 x_6^4 \in \leadmon & x_6^4 & \{ 1 \}\vspace{.05in}\\
         x_6  & x_6^5 \in \leadmon & x_6^4 & \{ 1 \}
    \end{array}.
\]
}
By Corollary \ref{cor:Z''Cor} the image of the map $\fn{\EuScript{F}}{U}{\Hilb^{32}_{\mathbb{A}^6_{\gf}}}$ satisfies
\[
    \dim(\EuScript{F}(U)) \geq |\specl| + |\mathcal{Z}''| = (\lambda \cdot \tau) + (11 \cdot 3 + 6 +2)  = 17 \cdot 10 + 41 = 211.
\]
It now follows from Proposition \ref{prop:genericCond} that $[I]$ is a smooth point on an elementary component of dimension $211$, so that $I$ is a generic ideal.


\section{\SLICap\ ideals of shape $(n,\kappa,r,s)$} \label{sec:idealsOfShapeNkRs}

Our last main goal, accomplished in Section \ref{sec:PlausArgs}, is to develop a numerical criterion for picking out shapes $(n,\kappa,r,s)$ for which sufficiently general \SLI\ ideals of that shape are likely to be generic.   We make preparations in this section by discussing such \SLI\ ideals in detail.  We assume $n \geq 3$, $1 < \kappa < n$, and $2 \leq r < s$.

\subsection{The sets $\mathcal{O}$, $\leadmon$, and $\trailmon$} \label{ssec:FamOfExamplesBasics} The order ideal $\mathcal{O}$ is the lex-segment complement given by $\mathcal{O}$ $=$ $\cup_{d=0}^{s} \mathcal{O}_d$ (recall that $\mathcal{O}_{d}$ denotes the set of basis monomials of degree $d$), where  
\[
    \mathcal{O}_d\ = \ \left\{ 
    \begin{array}{l}
        \{ \text{monom's of degree } d \text{ in } x_1, \dots, x_n \},\ \text{if } 0 \leq d < r,\vspace{.05in}\\
        \{ \text{monom's of degree } d \text{ in } x_{n-\kappa + 1}, \dots, x_n\},\ \text{if } r \leq d \leq s,\vspace{.05in}\\
        \emptyset, \text{ if } d > s.
    \end{array} \right.
\]
We call $x_1$, \dots, $x_{n-\kappa}$ the \textbf{front variables}, and $x_{n-\kappa+1}, \dots, x_n$ the \textbf{back variables}; we then define the \textbf{front degree} (resp.\ \textbf{back degree}) of a monomial $m$ to be the sum of the exponents to which the front (resp.\ back) variables appear in $m$.  The sets of leading and trailing monomials are selected as follows:
\[
    \begin{array}{rcl}
        \leadmon & = & \{ \text{monom's in } x_1, \dots, x_n \text{ of deg. } r \} \setminus \mathcal{O}_d,\vspace{.05in}\\
              {} & = & \{x_1^r,\, x_1^{r-1}x_2,\,\dots, x_{n-\kappa}\, x_n^{r-1} \} = \{b_1,\, b_2,\, \dots,\, b_{\lambda}   \},\text{ and}\vspace{.05in}\\
        \trailmon & = & \mathcal{O}_s  = \{x_{n-\kappa + 1}^s,\, x_{n-\kappa +1}^{s-1}\, x_{n-\kappa + 2},\, \dots,\, x_n^s \}.
    \end{array}
\] 
It follows easily that
\begin{equation} \label{eqn:lambdaAndTau}
  \begin{array}{rcl}
    \lambda & = & |\leadmon| = {n-1+r \choose r} - {\kappa - 1 + r \choose r},\vspace{.05in}\\
    \tau    & = & |\trailmon| = {\kappa - 1 + s \choose s}, \text{ and}\vspace{.05in}\\
    \mu     & = & \sum_{d=0}^{r-1} {n-1+d \choose d} + \sum_{d=r}^s {\kappa-1 + d \choose d}.
  \end{array}
\end{equation}
Recall that a \SLI\ ideal $I$ built using these sets $\mathcal{O}$, $\leadmon$, and $\trailmon$ is said to have \textbf{shape} $(n,\kappa,r,s)$.  Examples \ref{ssec:(1,5,3,4)}, \ref{ssec:(1,5,3,4,5,6)}, \ref{ssec:(1,6,6,10)}, and \ref{ssec:(1,6,21,10,15)} are all of this form, having shapes $(5,2,2,3)$, $(5,2,2,5)$, $(6,3,2,3)$, and $(6,3,3,4)$, respectively.


\subsection{Boundary monomials} \label{ssec:plausArgsBdryMonoms}

Recalling that $\partial \mathcal{O}_d$ $\subseteq$ $\partial \mathcal{O}$ denotes the subset of degree-$d$ boundary monomials, one sees easily that $\partial \mathcal{O} = \cup_{d=r}^{s+1} \partial \mathcal{O}_d$, where 
\begin{equation} \label{eqn:listOfBdryMons}
    \begin{array}{rcl}
        \partial \mathcal{O}_{r} & = &  
\leadmon\ = \ \{b_1, b_2, \dots, b_{\lambda} \},\vspace{.05in}\\    
        \partial \mathcal{O}_d  & = & \{x_{\alpha} \cdot t_i \mid 1 \leq \alpha \leq n-\kappa,\ t_i \in \mathcal{O}_{d-1} \},\ r+1 \leq d \leq s, \text{ and}\vspace{.05in}\\
        \partial \mathcal{O}_{s+1} & = &  
             \{x_{\alpha} \cdot t_i \mid 1 \leq \alpha \leq n-\kappa,\ t_i \in \mathcal{O}_{s} \}\ \cup \vspace{.05in}\\ 
    {} & {} & \left \{ 
            \text{monomials of degree } s+1 \text{ in } x_{n-\kappa+1}, \dots, x_n  
              \right\} \vspace{.05in}\\
             {}                 & = & \partial \trailmon .
    \end{array} 
\end{equation}
Consequently, $|\partial \mathcal{O}|$ $=$ $\nu$ is given by
{\Small
\[
    \nu\ =\ {n -1 + r \choose r} - {\kappa-1+r \choose r} + \left( \sum_{d=r+1}^{s+1}  (n-\kappa){\kappa-1+d-1 \choose d-1} \right) + {\kappa + s \choose s+1} .
\]
}
\begin{rem} \label{rem:Cond1ForEffIsOK}
    Condition (i) of Proposition {\rm \ref{prop:effTestingProp}}, namely, that every non-leading boundary monomial $b_j$ is a multiple of a monomial in $Q$ $=$ $\partial \leadmon\, \cup\, \partial \trailmon$, is easily seen to hold in this case.  Indeed, 
\[
    \partial \leadmon = \{ \text{monom's of degree } d+1 \text{ and front degree} \geq 1 \},
\]
and clearly every non-leading boundary monomial with front degree $\geq$ $1$ is divisible by a monomial in $\partial \leadmon$.  The only other non-leading boundary monomials are the monomials of degree $s+1$ and front degree $0$, and these all lie in $\partial \trailmon$.
\end{rem}

\subsection{Linear syzygies} \label{ssec:neighborSyz}

In this section we compute the cardinality of the set $\tarmons$ of target monomials (\ref{eqn:tarMonDef}), from which, by (\ref{eqn:linSyzCount}), we obtain the dimension of the $\gf$-vector space of linear syzygies $\psi$ $=$ $(n+1) \cdot \nu - |\tarmons|$.   Writing $\tarmons_d$ $\subseteq$ $\tarmons$ for the subset of monomials of degree $d$, we list the elements of 
$\tarmons = \cup_{d=r}^{s+2}\, \tarmons_d$ degree-by-degree: 
First, it is clear that 
\[
    \tarmons_r =  \partial \mathcal{O}_r, \text{ so }\ 
  |\tarmons_r| = \lambda = {n-1+r \choose r} - {\kappa - 1 + r \choose r}.
\]
Next, since 
\[
    \cup_{\alpha = 1}^{n} (x_{\alpha}\cdot \partial \mathcal{O}_r) \subseteq \tarmons_{r+1} \subseteq \{\text{monom's of degree } r+1 \text{ and front deg.} \geq 1   \},
\]
and the extremes are clearly the same, we have that
\begin{equation} \label{eqn:tarmonsRplus1}
        |\tarmons_{r+1}| = {n + r \choose r+1} - {k + r \choose r+1}.
\end{equation}
One checks easily that, for $r+1$ $\leq$ $d$ $\leq$ $s$, 
\[
     \partial \mathcal{O}_{d} = \{ \text{monom's of degree } d\text{ and front degree }1 \}.
\]
Hence, for $r+2$ $\leq$ $d$ $\leq$ $s$, we have that
\[
  \begin{array}{rcl}
    \tarmons_d & = & \left( \partial \mathcal{O}_{d}\right) \cup \left( \cup_{\alpha=1}^{n}\, x_{\alpha}\cdot \partial \mathcal{O}_{d-1}\right)\vspace{.05in}\\
   {} & =  &\{m \mid \deg(m) =d\ \& \operatorname{front-deg}(m) = 1 \text{ or }2  \}, \text{ so}\vspace{.05in}\\
 |\tarmons_d| & = & (n-\kappa)\cdot {\kappa + d - 2 \choose d-1} + {n-\kappa + 1 \choose 2}\cdot {\kappa + d -3 \choose d-2}.
  \end{array}
\] 
Continuing, we next observe that
\[ 
    \begin{array}{rcl}
       \tarmons_{s+1} & = & \partial \mathcal{O}_{s+1} \cup \left( \cup_{\alpha=1}^{n}\, x_{\alpha}\cdot \partial \mathcal{O}_{s}\right)\vspace{.05in}\\  
     {} &  = & \left\{
                  \begin{array}{c}\text{monom's of deg. } s+1 \text{ and}\vspace{.05in}\\
\text{front deg. }0,\, 1, \text{ or } 2
                  \end{array} \right\}, \text{ so}\vspace{.05in}\\
       |\tarmons_{s+1}| & = & {\kappa  + s \choose s+1} + (n-\kappa)\cdot {\kappa-1+s \choose s} + {n-\kappa + 1 \choose 2}\cdot {\kappa-2 +s \choose s-1}.
    \end{array}
\]
Finally, we observe that 
\[
    \begin{array}{c}
        \tarmons_{s+2} = \cup_{\alpha=1}^{n}\, x_{\alpha}\cdot \partial{\mathcal{O}}_{s+1} = \left\{
          \begin{array}{c}
            \text{monom's of deg. } s+2\text{ and }\vspace{.05in}\\
            \text{front deg. } 0,\, 1, \text{ or } 2
         \end{array} \right \}, \text{ so }\vspace{.05in}\\
        |\tarmons_{s+2}|  = {\kappa  + s+1 \choose s+2} + (n-\kappa)\cdot {\kappa + s \choose s+1} + {n-\kappa + 1 \choose 2}\cdot {\kappa-1 +s \choose s},
    \end{array}
\]
and $|\tarmons|$ $=$ $\sum_{d=r}^{s+2}|\tarmons_d|$.  By inspection of the preceding results, we obtain the following
\begin{lem} \label{lem:asymptotics}
    Let $I$ be a \SLI\ ideal of shape $(n,\kappa,r,s)$.  If we hold $\kappa$, $r \geq 2$, and $s$ $>$ $r$ constant, and allow $n$ to increase, then the quantities $\mu$, $\lambda$, $\tau$, $\nu$, $|\tarmons|$, and $\psi$ are polynomials in $n$ with the following dominant terms:
\[
    \begin{array}{rcl}
        \mu & \approx & \frac{n^{r-1}}{(r-1)!},\vspace{.05in}\\
        \lambda & \approx & \frac{n^r}{r!},\vspace{.05in}\\
        \tau & = & {\kappa-1+s \choose s}\ =\  \text{\rm constant in } n,\vspace{.05in}\\
        \nu & \approx & \frac{n^r}{r!},\vspace{.05in}\\
        |\tarmons| & \approx & \frac{n^{r+1}}{(r+1)!},\vspace{.05in}\\
        \psi      & \approx & n\cdot \frac{n^r}{r!} -\frac{n^{r+1}}{(r+1)!} = \frac{r}{(r+1)!}\cdot n^{r+1} .
    \end{array}
\] \qed
\end{lem}


\subsection{A linearly independent set in $\tansp_{[I]}$} \label{ssec:LinIndepTanVecSetInMainCase}

Since the order ideal $\mathcal{O}$ under consideration is a lex-segment complement, we choose the sets $\Delta'_{\alpha}$ as described in Section \ref{ssec:ZDiscussion}.  One checks easily that, for $1 \leq \alpha \leq n-\kappa$, 
\[
    \begin{array}{l}
         b_{j_{\alpha}} = x_{\alpha}\,x_n^{r-1} \in \leadmon,\ t_{i_{\alpha}} = x_n^{r-1}, \text{ and }\vspace{.05in}\\ 
 \Delta'_{\alpha} = \{\text{monom's } m \in \gf[x_{n-\kappa+1}, \dots, x_n] \mid 0 \leq \deg(m) \leq s-r\}, 
    \end{array}     
\]
and, for $n-\kappa+1 \leq \alpha \leq n$,
\[
         b_{j_{\alpha}} = x_{\alpha}\,x_n^{s} \in \partial \mathcal{O} \setminus \leadmon,\ t_{i_{\alpha}} = x_n^s \in \trailmon, \text{ and }\  \Delta'_{\alpha} = \{ 1 \}.
\]

\begin{lem} \label{lem:SCupZ'LinIndInShapeCases}
    For any \SLI\ ideal $I = I_p$ of shape $(n,\kappa,r,s)$  the set $\EuScript{S}_p \cup \EuScript{Z'}_p$ $\subseteq$ $\tansp_{[I]}$ is linearly independent. Consequently, if $\dim_{\gf}(\tansp_{[I]})  =   |\specl| + |\mathcal{Z}'|$  
(equivalently, if $\EuScript{S}_p \cup \EuScript{Z}'_p$ is a $\gf$-basis of $\tansp_{[I]}$), then $I$ is generic.
\end{lem}

\proof
    The first statement results from Corollary \ref{cor:GenZPropCor2}.  The second statement then follows from Proposition \ref{prop:genericCond}.
\qed
\medskip

Henceforth we say that the \SLI\ ideal $I$ of shape $(n,\kappa,r,s)$ is 
\textbf{shape-generic}
if and only if the following equation holds:
{\Small
\begin{equation} \label{eqn:tanSpDimInGenShapeCase}
    \begin{array}{c}
      \dim_{\gf}(\tansp_{[I]})\ =\ |\specl|\ +\ |\mathcal{Z}'|\ =\  \lambda \cdot \tau + (n-\kappa) \cdot |\Delta'_1| + \kappa \cdot |\Delta'_{n-\kappa+1}|\ = \vspace{.05in}\\
      \left( {n-1+r \choose r} - {\kappa - 1\ +\ r \choose r} \right) \cdot {\kappa - 1 + s \choose s}\ +\ (n-\kappa) \cdot \left( \sum_{d=0}^{s-r}{\kappa-1+d \choose d}  \right)\ +\ (\kappa)\cdot 1
    \end{array} .
\end{equation}
}
\subsection{Tangent space relations} \label{ssec:tanSpRelnCt}

Recall from Section \ref{sec:tanSpace} that a tangent vector at $[I]$ is given by an $R$-homomorphism $\fn{v}{I}{R/I}$, which is determined by the images 
\[
    v(g_j)\ =\ \sum_{i=1}^{\mu}a_{ij}t_i \in \Span_{\gf}(\mathcal{O}),\ 1 \leq j \leq \nu,
\]
where the $g_j$ are the elements of the $\mathcal{O}$-border basis of $I$.  Given a syzygy $\sum_{j=1}^{\nu}f_j g_j$ $=$ $0$, we obtain the linear relations on the $a_{ij}$ as described in Section \ref{ssec:tanSpGen}:
\[
  \begin{array}{rcl}
    0\ =\ v \left( \sum_{j=1}^{\nu}f_j\, g_j \right)\ =\ \sum_{j=1}^{\nu}f_j\, v(g_j) &  = & \sum_{j=1}^{\nu}f_j\cdot(\sum_{i=1}^{\mu}a_{ij} t_i)\vspace{.05in}\\
{} &  \equiv & \sum_{i=1}^{\mu}\mathbf{b}^{(f_j)}_{i}t_{i} \text{ mod } I,
  \end{array}
\]
where each of the coefficients $\mathbf{b}^{(f_j)}_{i}$ is a $\gf$-linear combination of the $a_{ij}$ that must vanish.  Viewing the $a_{ij}$ as indeterminates, we represent a ``generic tangent vector'' as the $(\mu \nu)$-tuple $(a_{1,1},\, a_{2,1},\, \dots,\, a_{\mu \lambda})$, as in (\ref{eqn:TanVec}).  By computing the tangent space relations $\mathbf{b}_{i}^{(f_j)}$ for a basis of the linear syzygies, and then computing the dimension of the vector space that they span, we obtain  $\dim_{\gf}(\tansp_{[I]})$ as in (\ref{eqn:tanSpDimComp}).    

We assign a \textbf{degree} to each indeterminate $a_{ij}$ and each relation $\mathbf{b}_{i}^{(f_j)}$ as follows: 
\[
    \deg(a_{ij})\ =\ \deg(\mathbf{b}_{i}^{(f_j)})\ =\ \deg(t_{i}).
\]  

Our goal here is to identify and count the $a_{i'j'}$ that can appear non-trivially in the relation $\mathbf{b}_{i}^{(f_j)}$ associated to a linear syzygy $(f_j)$. To this end, let $f_{j',k}$ denote one of the terms of $f_{j'}$, so $f_{j',k}$ is either a constant or a constant times a variable, and observe that the product of a variable $x_{\alpha}$ with a basis monomial $t_{i'}$ of degree $d$ is a monomial $m$ $=$ $x_{\alpha}\, t_{i'}$ of degree $d+1$ such that exactly one of the following holds:
\begin{equation} \label{eqn:CongruencePossibilities}
    \begin{array}{rcl}
        m \in \mathcal{O}_{d+1} & \Rightarrow & m \equiv m\ (\text{mod } I),\vspace{.05in}\\ 
        m \in \leadmon & \Rightarrow & m \equiv N \in \Span_{\gf}(\trailmon)\ (\text{mod } I),\vspace{.05in}\\
        m \in \partial \mathcal{O} \setminus \leadmon & \Rightarrow & m \equiv 0\ (\text{mod } I).
    \end{array} 
\end{equation}
We see that the indeterminate $a_{i'j'}$ of degree $d$ can appear in the relation $\mathbf{b}_i^{(f_j)}$ if one of the following holds (here $0 \neq c \in \gf$):
\begin{itemize}
    \item $f_{j',k} = c $, and $i' = i$, so that $\deg(\mathbf{b}_i^{(f_j)}) = d$;
    \item $f_{j',k} = c x_{\beta}$ and $x_{\beta}\cdot t_{i'} = t_i$ $\in$ $\mathcal{O}$,  so that  $\deg(\mathbf{b}_i^{(f_j)})  = d+1$;
    \item $f_{j',k} = c x_{\beta}$ and $x_{\beta} \cdot t_{i'} \in \leadmon$, so that $d = r-1$ and $\deg(\mathbf{b}_i^{(f_j)})  = s$. 
\end{itemize} 

\begin{rem} \label{rem:syzCoeffsRem}

	Let $I = \{g_j \}$ be a \SLI\ ideal of shape $(n,\kappa,r,s)$.  For $1 \leq j' \leq \lambda$, a term $f_{j',k}$ in the $j'$-entry of a linear syzygy $(f_j)$ can never equal a non-zero constant $c$, since $c \cdot b_{j'}$ can never cancel out of the expression $\sum_{j=1}^{\nu}f_j\cdot g_j$ for $b_{j'}$ $\in$ $\leadmon$.    

\end{rem}

For each degree $0$ $\leq$ $d$ $\leq$ $s$, we let $A_d$ denote the set of indeterminates $a_{ij}$ that appear in at least one of the tangent space relations of degree $d$. We proceed to identify the possible members of these sets and find upper bounds on their cardinalities, for which it is convenient to adopt the following terminology: we say that $a_{ij}$ $\in$ $A_d$ \textbf{stays put} if $\deg(a_{ij})$ $=$ $d$, \textbf{moves up by 1} if $\deg(a_{ij})$ $=$ $d-1$, or \textbf{jumps up} if $\deg(a_{ij})$ $=$ $r-1$ and $d=s$; in addition, we take $|\mathcal{O}_{-1}|$ $=$ $0$.  Summarizing the foregoing observations, we obtain the following

\begin{lem} \label{lem:ACountLem}
    The $a_{ij}$ $\in$ $A_d$ that stay put have indices satisfying $b_j$ $\in$ $\partial \mathcal{O}\setminus \leadmon$ and $t_i$ $\in$ $\mathcal{O}_d$, for $0 \leq d \leq s$.  The $a_{ij}$ $\in$ $A_d$ that move up by $1$ have indices satisfying $b_j$ $\in$ $\partial \mathcal{O}$ and $t_i$ $\in$ $\mathcal{O}_{d-1}$, for $0$ $\leq$ $d$ $\leq$ $r-1$ and $r+1$ $\leq$ $d$ $\leq$ $s$, and, for $d = r$, indices satsfying $b_j$ $\in$ $\partial \mathcal{O}$ and $t_i$ $\in$ $\{m \in \mathcal{O}_{r-1} \mid \operatorname{back-deg.}(m) = r-1 \}$.  Finally, the $a_{ij}$ $\in$ $A_s$ that jump up have indices satisfying $b_j$ $\in$ $\partial \mathcal{O}$ and $t_i$ $\in$ $\mathcal{O}_{r-1}$.  Whence,
\[
        |A_d| \ \leq \ \left\{
                       \begin{array}{l}
                            (\nu - \lambda) \cdot {n-1+d \choose d} + \nu \cdot {n-1+d-1 \choose d-1}, \text{ if } 0 \leq d \leq r-1,\vspace{.05in}\\
                            (\nu-\lambda) \cdot {\kappa - 1 + r \choose r} + \nu \cdot {\kappa-1+r-1 \choose r-1 }, \text{ if } d = r,\vspace{.05in}\\
                            (\nu - \lambda) \cdot {\kappa - 1 + d \choose d} + \nu \cdot {\kappa - 1 + d-1 \choose d-1}, \text{ if } r+1 \leq d \leq s-1,\vspace{.05in}\\
                            (\nu - \lambda) \cdot {\kappa - 1 + s \choose s} + \nu \cdot {\kappa - 1 + s-1 \choose s-1} + \nu \cdot {n-1+ r-1 \choose r-1}, \text{ if }  d = s.
                       \end{array} \right.
\]
\qed
\end{lem}

\subsection{Quasi-efficiency} \label{ssec:quasi-eff}

By Remark \ref{rem:GenButNotEff}, a shape-generic ideal $I$ of shape $(n,\kappa,r,s)$ need not be efficient.  However, as we show in this section, it must have the following property that we call \textbf{quasi-efficiency}: Every tangent vector $\fn{v}{I}{R/I}$ in $\tansp_{[I]}$ is determined by the images $v(g_j)$ for $1\leq j \leq \lambda$.  Clearly efficiency implies quasi-efficiency, since if $I$ is generated by $g_1, \dots, g_{\lambda}$, then any $R$-homomorphism of $I$ is determined by the images of these generators. 

\begin{lem} \label{lem:degOfALowerBd} Let $v$ $\in$ $\EuScript{S}_p \cup \EuScript{Z}'_p$, and let $(a_{ij})$ denote the associated tuple.  Then the minimal degree of a non-zero component $a_{ij}$ is as follows:
\begin{itemize}
  \item If $v$ $=$ $v_{p,ij}$ $\in$ $\EuScript{S}_p$, then $v$ has a single non-zero component of degree $s$.  
  \item If $v$ $=$ $v_{p,\alpha,\delta}$ $\in$ $\EuScript{Z}'_p$, then the minimal degree of a non-zero component in $(a_{ij})$ is given by 
      \begin{itemize}
        \item[$\ast$] $\deg(t_{i_{\alpha,\delta},j_{\alpha}})  =  r-1 + \deg(m_{\alpha,\delta}) \leq s-1, \text{ if } x_{\alpha}$ is a front variable, and
        \item[$\ast$] $r-1, \text{ if } x_{\alpha}$ is a back variable. 
        \end{itemize} 
\end{itemize}
Moreover, when $x_{\alpha}$ is a back variable, the components of $v_{p, \alpha, \delta}$ can only be non-zero in degrees $r-1$, $s-1$,  and $s$.
\end{lem} 

\proof
     The first bulleted statement is immediate from Section \ref{ssec:TVFam1}: $v_{p,ij}$ $\in$ $\EuScript{S}_p$ has a single non-zero entry with \SLI\ index pair $ij$, which implies that $t_i$ $\in$ $\trailmon$.  In our context, this yields $\deg(a_{ij})$ $=$ $\deg(t_i)$ $=$ $s$.   

We now prove the second bulleted statement assuming that $x_{\alpha}$ is a front variable.  From the proof of Proposition \ref{prop:GeneralZProp}, we have that the $(i_{\alpha,\delta},j_{\alpha})$-component of $v$ $=$ $v_{p,\alpha,\delta}$ is non-zero --- indeed, it is the only non-zero component with index of the form $(i,j_{\alpha})$ --- and this component has degree $\deg(t_{i_{\alpha,\delta}})$.  So it remains to show that a non-zero component $a_{i,j}$ with $j$ $\neq$ $j_{\alpha}$ cannot have a strictly smaller degree.   

    Recall that, by (\ref{eqn:ZTanVecHomo}) the tangent vector $\fn{v}{I}{R/I}$ is given by 
\[
     g_j\ \mapsto \frac{\partial g_j}{\partial x_{\alpha}} \cdot m_{\alpha,\delta}\ \equiv \ \sum_{i=1}^{\mu}a_{ij}t_i\ (\text{mod } I),\ \ 1 \leq j \leq \nu.
\]
Since $x_{\alpha}$ is a front variable and the trailing monomials only involve back variables, we have that $\frac{\partial g_j}{\partial x_{\alpha}}$ $=$ $\frac{\partial b_j}{\partial x_{\alpha}}$, so $\frac{\partial g_j}{\partial x_{\alpha}} \cdot m_{\alpha,\delta}$ is either equal to $0$ or to $\text{(const)}\cdot M$, where the monomial $M$ $=$ $\frac{b_{j}}{x_{\alpha}} \cdot m_{\alpha,\delta}$ has degree $\deg(b_j) - 1 + \deg(m_{\alpha,\delta})$ $\geq$ $r-1+\deg(m_{\alpha,\delta})$.  Modulo $I$, we have that $M$ is congruent to one of the following: 
\begin{equation} \label{eqn:3Possibilities}
    \begin{array}{l}
        M, \text{ if } M \in \mathcal{O},\vspace{.05in}\\
        N \in \Span_{\gf}(\trailmon)\   \text{ (of degree $s$), if } M \in \leadmon,\text{ or }\vspace{.05in}\\
        0, \text{ if } M \notin \mathcal{O} \text{ and } M \notin \leadmon.
    \end{array}
\end{equation} 
In each case, we see that no non-zero component of degree $<$ $r-1+ \deg(m_{\alpha,\delta})$ occurs.  

Now let $x_{\alpha}$ be a back variable (which implies that $m_{\alpha,\delta}$ $=$ $m_{\alpha,1}$ $=$ $1$), and consider the components coming from $\frac{\partial g_j}{\partial x_{\alpha}} \cdot 1$.  If $b_j$ $\in$ $\leadmon$, then non-zero components of degree $s-1$ can arise from $-\frac{\partial N_j}{\partial x_{\alpha}}$, where $g_j\ = \ b_j - N_j$.  Otherwise (in all cases), non-zero components can result from $\frac{\partial b_j}{\partial x_{\alpha}}$ $=$ $\text{(const)}\cdot M$ (if non-zero).  The three possibilities (\ref{eqn:3Possibilities}) again present themselves; however, because $x_{\alpha}$ is a back variable, we can say more.  Consider the boundary monomials $b_j$ such that $\frac{b_j}{x_{\alpha}}$ $=$ $M$ $\in$ $\mathcal{O}$.  A moment's reflection shows that this is possible only in one of the following two cases, both of which occur: $\deg(b_j)$ $=$ $r$ and is divisible by $x_{\alpha}$, or $\deg(b_j)$ $=$ $s+1$ and $b_j$ involves only back variables (including $x_{\alpha}$), since for all other degrees $r+1$ $\leq$ $d$ $\leq$ $s$, $b_j$ involves a front variable, hence so does $M$ (of degree $\geq r$) $\Rightarrow$ $M$ $\notin$ $\mathcal{O}$.  It follows that the only possible degrees for non-zero components are $r-1$, $s-1$, and $s$, as the final statement asserts, and the minimum $r-1$ is attained. This completes the proof of the lemma.
\qed

\begin{cor} \label{cor:AsAre0InDegsLEQN-2}
    Let $I=I_p$ be shape-generic, and let $v$ $=$ $(a_{ij})$ $\in$ $\tansp_{[I]}$. Then $\deg(a_{ij})$ $\leq$ $r-2$ $\Rightarrow$ $a_{ij}$ $=$ $0$.
\end{cor}

\proof By definition, if $I$ $=$ $I_p$ is shape-generic, then the set of vectors $\EuScript{S}_p \cup \EuScript{Z}'_p$ is a basis of $\tansp_{[I_p]}$.  Proposition \ref{lem:degOfALowerBd} implies that the desired conclusion holds for this basis; whence, it holds for all $v \in \tansp_{[I]}$.  
\qed

\begin{rem} \label{rem:SimToIarrobEmsalemResult}
    The preceding Corollary extends a consequence of \cite[Lemma 2.31, p.\ 162]{Iarrob-Emsalem} to our case.
\end{rem}

Recall from Corollary \ref{cor:GenZPropCor1A} that the $(i_{\alpha,\delta},j_{\alpha})$-component of $v_{p,\alpha,\delta}$ is non-zero, and for all $\beta$ $\leq$ $\alpha$ and $(\beta,\delta')$ $\neq$ $(\alpha,\delta)$, the $(i_{\alpha,\delta},j_{\alpha})$-component of $v_{p,\beta,\delta'}$ $=$ $0$. In this sense, the $(i_{\alpha,\delta},j_{\alpha})$-component of $v_{p,\alpha,\delta}$ acts as a ``characteristic function'' for  $v_{p,\beta,\delta'}$.
  When $x_{\alpha}$ is a front variable, $b_{j_{\alpha}}$ $=$ $x_{\alpha}x_n^{r-1}$ $\in$ $\leadmon$, but when $x_{\alpha}$ is a back variable, $b_{j_{\alpha}}$ $=$  $x_{\alpha}x_n^{s}$ $\in$ $\partial \trailmon$.  To prove quasi-efficiency, we use the following modified set of ``characteristic function'' components of index  $(\hat{i}_{\alpha,\delta},\hat{j}_{\alpha})$ such that $b_{\hat{j}_{\alpha}}$ $\in$ $\leadmon$  for all $1$ $\leq$ $\alpha$ $\leq$ $n$:
\[
    \begin{array}{l}
        b_{\hat{j}_{\alpha}} = b_{j_{\alpha}},\ t_{\hat{i}_{\alpha,\delta}} = t_{i_{\alpha,\delta}} = x_n^{r-1}\cdot m_{\alpha,\delta}, \text{ if } x_{\alpha} \text{ is a front variable, and}\vspace{.05in}\\
        b_{\hat{j}_{\alpha}} = x_{\alpha}x_1^{r-1},\ t_{\hat{i}_{\alpha,\delta}} = t_{\hat{i}_{\alpha,1}} = x_1^{r-1}\cdot 1, \text{ if } x_{\alpha} \text{ is a back variable}.
    \end{array}
\]  

\begin{lem} \label{lem:quasiEffPrepLem}
    Let $I = I_p$ be a \SLI\ ideal of shape $(n,\kappa,r,s)$.  Then for all variables $x_{\alpha}$ and all  $1$ $\leq$ $\delta$ $\leq$ $|\Delta'_{\alpha}|$, we have the following:
\begin{itemize}
  \item[i.] The $(\hat{i}_{\alpha,\delta},\hat{j}_{\alpha})$-component of $v_{p,\alpha,\delta}$ is non-zero.
  \item[ii.] For all variables $x_{\beta}$ and all $1$ $\leq$ $\delta'$ $\leq$ $\Delta'_{\beta}$, if $\beta$ $\leq$ $\alpha$ and $(\beta,\delta')$ $\neq$ $(\alpha,\delta)$, then the $(\hat{i}_{\alpha,\delta},\hat{j}_{\alpha})$-component of $v_{p,\beta,\delta'}$ $=$ $0$.
\end{itemize} 
\end{lem}

\proof
     The truth of the first statement for the front variables $x_{\alpha}$ was shown in the proof of Proposition \ref{prop:GeneralZProp}, so suppose that $x_{\alpha}$ is a back variable.  By (\ref{eqn:ZTanVecHomo}), the $(\hat{i}_{\alpha,1},\hat{j}_{\alpha})$-component of $v_{p,\alpha,1}$ is the coefficient of $t_{\hat{i}_{\alpha,1}}$ $=$ $x_1^{r-1}$ in $\frac{\partial g_{\hat{j}_{\alpha}}}{\partial x_{\alpha}} \cdot 1$ modulo $I$.  But $g_{\hat{j}_{\alpha}}$ $=$ $x_{\alpha} x_1^{r-1} - N_{\hat{j}_{\alpha}}$, where $N_{\hat{j}_{\alpha}}$ $\subseteq$ $\Span_{\gf}(\trailmon)$, so it is clear that the desired coefficient is $1$.

To prove the second statement, we must show, given $\beta$ $\leq$ $\alpha$ and $(\beta,\delta')$ $\neq$ $(\alpha,\delta)$, that the coefficient of $t_{\hat{i}_{\alpha,\delta}}$ in $\frac{\partial g_{\hat{j}_{\alpha}}}{\partial x_{\beta}}\cdot \delta'$ modulo $I$ is $0$.  We first consider the case in which $x_{\alpha}$ is a back variable, so that $\Delta'_{\alpha}$ $=$ $\{ 1 \}$, and let $(x_{\beta},\delta')$ $\neq$ $(x_{\alpha},1)$, which implies that $\beta$ $<$ $\alpha$.  Then
\[
    \frac{\partial g_{\hat{j}_{\alpha}}}{\partial x_{\beta}} \cdot \delta'\ =\ \frac{\partial (x_{\alpha} x_1^{r-1})}{\partial x_{\beta}}\cdot \delta' - \frac{\partial N_{\hat{j}_{\alpha}}}{\partial x_{\beta}} \cdot \delta'.
\] 
The first term on the RHS of the preceding equation is $0$ provided $x_{\beta}$ $\neq$ $x_1$, and if $x_{\beta}$ $=$ $x_1$, it equals $(r-1) x_{\alpha} x_1^{r-2} \cdot \delta'$ $=$ $(r-1)\cdot m$, where $m$ is a monomial of degree $\geq r-1$.  There are three possibilities for the value of $m$ modulo $I$, as in (\ref{eqn:CongruencePossibilities}), and a moment's reflection shows that none of these possibilities can include a non-zero multiple of $t_{\hat{i}_{\alpha,1}}$ $=$ $x_1^{r-1}$.  Similarly, the second term on the RHS consists of a linear combination of monomials of degree $\geq s-1$, and we again conclude that, modulo $I$, no non-zero multiple of $x_1^{r-1}$ can appear.  It follows that statement ii.\ holds when $x_{\alpha}$ is a back variable.

Now consider the case in which $x_{\alpha}$ is a front variable. One can then obtain the desired conclusion immediately from Corollary \ref{cor:GenZPropCor1A} or from the following argument: 
Choose $\beta$ $\leq$ $\alpha$ (so $x_{\beta}$ is a front variable) and $(\beta,\delta')$ $\neq$ $(\alpha,\delta)$, and compute the coefficient of $t_{\hat{i}_{\alpha,\delta}}$ $=$ $x_n^{r-1}\cdot m_{\alpha,\delta}$ in
\[
  \begin{array}{rcl}
    \frac{\partial g_{\hat{j}_{\alpha}}}{\partial x_{\beta}} \cdot m_{\beta,\delta'} & = & \left( \frac{\partial b_{\hat{j}_{\alpha}}}{\partial x_{\beta}} - \frac{\partial N_{\hat{j}_{\alpha}}}{\partial x_{\beta}} \right) \cdot m_{\beta,\delta'}\vspace{.05in}\\
{} & = & \left( \frac{\partial\, x_{\alpha}\, x_{n}^{r-1}}{\partial x_{\beta}} - \frac{\partial N_{\hat{j}_{\alpha}}}{\partial x_{\beta}}\right) \cdot m_{\beta,\delta'} \text{ mod } I.
  \end{array}
\] 
If $\beta$ $<$ $\alpha$, the last expression is clearly $0$, and if $\beta$ $=$ $\alpha$, it equals $t_{\hat{i}_{\alpha,\delta'}}$.  This can yield a non-zero coefficient for $t_{\hat{i}_{\alpha,\delta}}$ only if $\delta$ $=$ $\delta'$, which is ruled out by the hypothesis that $(\beta,\delta')$ $\neq$ $(\alpha,\delta)$.  This completes the proof.
\qed
\medskip

We are now ready to prove that a shape-generic \SLI\ ideal $I$ is quasi-efficient.

\begin{prop} \label{prop:gsDetermineTanVec}
    Let $I = I_p$ be a shape-generic \SLI\ ideal of shape $(n,\kappa,r,s)$, and let $\fn{v}{I}{R/I}$ be a tangent vector at $[I]$ with associated tuple $(a_{ij})$.  Then $v$ is determined by the images $v(g_j)$ for $1\leq j \leq \lambda$. 
\end{prop}

\proof
     We begin by writing $v$ as a (unique) linear combination of the elements of the basis $\EuScript{S}_p$ $\cup$ $\EuScript{Z'}_p$:
\begin{equation} \label{eqn:basisExp}
    v\ =\ \sum_{v_{p,ij} \in \EuScript{S}_p} d_{ij}\,v_{p,ij} + \sum_{v_{p,\alpha,\delta} \in \EuScript{Z'}_p} d_{\alpha,\delta}\, v_{p,\alpha,\delta},\ \  \ d_{ij},\ d_{\alpha,\delta} \in \gf . 
\end{equation}
     It suffices to show that the coefficients $d_{ij}$ and $d_{\alpha,\delta}$ are completely determined by $v(g_1)$, \dots, $v(g_{\lambda})$.  We begin by equating the $(\hat{i}_{n,1},\hat{j}_n)$-components on both sides of the equation.  By Lemma \ref{lem:quasiEffPrepLem}, we know that the $(\hat{i}_{n,1},\hat{j}_n)$-component of $v_{p,\beta,\delta'}$ is $0$ for all $\beta$ $\leq$ $\alpha$ and $(\beta,\delta')$ $\neq$ $(n,1)$, which includes all the pairs $(\beta,\delta')$ $\neq$ $(n,1)$.  Furthermore, the $(\hat{i}_{n,1},\hat{j}_n)$-components of the $v_{p,ij}$ are all $0$ since (by Lemma \ref{lem:degOfALowerBd}) the only non-zero component of $v_{p,ij}$ has degree $s$, and the degree of $t_{\hat{i}_{n,1}}$ $=$ $x_1^{r-1}$ is $r-1$ $<$ $s$. From this it follows that the coefficient $d_{n,1}$ is determined by the $(\hat{i}_{n,1},\hat{j}_n)$-component of $v$, which is the coefficient of $t_{\hat{i}_{n,1}}$ in $v(g_{\hat{j}_n})$.  This shows that $d_{n,1}$ is determined by   $v(g_1)$, \dots, $v(g_{\lambda})$.    

Proceeding by descending induction on $\alpha$, we assume that for some $1$ $\leq$ $\alpha$ $<$ $n$, all of the coefficients $d_{\beta,\delta'}$ for $\alpha+1$ $\leq$ $\beta$ $\leq$ $n$ are completely determined by  $v(g_1)$, \dots, $v(g_{\lambda})$ (and have been computed).  
We then equate the $(\hat{i}_{\alpha,\delta},\hat{j}_{\alpha})$-components on both sides of equation (\ref{eqn:basisExp}). 
Lemma \ref{lem:quasiEffPrepLem} implies that for all $\beta$ $\leq$ $\alpha$ and $(\beta,\delta')$ $\neq$ $(\alpha,\delta)$, the $(\hat{i}_{\alpha,\delta},\hat{j}_{\alpha})$-component of $v_{p,\beta,\delta'}$ is $0$, and the same is again true for all the $v_{p,ij}$, since none of the monomials $t_{\hat{i}_{\alpha,\delta}}$ can have degree $s$.  It follows that the value of $d_{\alpha,\delta}$ is determined by the coefficient of $t_{\hat{i}_{\alpha,\delta}}$ in $v(g_{\hat{j}_{\alpha}})$ and the previously-computed $d_{\beta,\delta'}$, so we conclude that for all $1$ $\leq$ $\delta$ $\leq$ $\Delta'_{\alpha}$, the coefficients $d_{\alpha,\delta}$ are determined by $v(g_1)$, \dots, $v(g_{\lambda})$.  It follows by induction that this is so for all the coefficients $d_{\alpha,\delta}$, $1 \leq \alpha \leq n$, $1 \leq \delta \leq |\Delta'_{\alpha}|$.

It is now clear that the remaining coefficients can be computed by equating the \SLI\ $(i,j)$-components on both sides of equation (\ref{eqn:basisExp}), so the value of each $d_{ij}$ can be computed from the coefficient of $t_{i,j}$ in $v(g_j)$  and the previously-computed values of the $d_{\alpha,\delta}$, hence is again determined by the values $v(g_1)$, \dots, $v(g_{\lambda})$, and we are done.
 \qed

\section{A criterion for plausible genericity} \label{sec:PlausArgs}

     To conclude the paper, we present a numerical criterion for identifying shapes $(n,\kappa,r,s)$ such that sufficiently general \SLI\ ideals associated to those shapes are likely to be shape-generic; we will call such shapes \textbf{plausible}.  

\subsection{The criterion} \label{ssec:PlausCrit}

Roughly speaking, the criterion is this: $(n,\kappa,r,s)$ is deemed plausible if the following two conditions hold:
\begin{description}
  \item[1] there are enough tangent space relations in each degree to allow the ranks of these sets of relations (if sufficiently general) to attain their maximum possible values, and
  \item[2] sufficiently general \SLI\ ideals $I$ of the given shape are likely to be $\vartheta$-efficient, and therefore likely to be quasi-efficient.
\end{description} 

We make these conditions computably precise and briefly argue for their
 reasonableness as follows: 

\begin{description}
  \item[1] Examples suggest that for shape-generic \SLI\ ideals $I$ of shape $(n,\kappa,r,s)$, the tangent space relations in each degree will attain (or nearly attain) their maximum possible ranks.  Of course, the rank of the tangent space relations in degree $d$ is bounded above by $|A_d|$, the number of indeterminates $a_{i,j}$ that appear in the relations of degree $d$, so we make condition 1 precise by requiring that the number of tangent space relations in each degree $0 \leq d \leq s$ is $\geq$ the upper bound on $|A_d|$ given in Lemma \ref{lem:ACountLem}.  Hence, if condition 1 holds, there are enough tangent space relations to render $I$ shape-generic, assuming that these relations are sufficiently independent.  
  \item[2] Since Proposition \ref{lem:quasiEffPrepLem} requires that any shape-generic \SLI\ ideal $I$ be quasi-efficient, and $\vartheta$-efficiency is an easy-to-check condition that implies quasi-efficiency, we require condition 2 in addition to condition 1.  In light of Remark \ref{rem:Cond1ForEffIsOK}, we know that \SLI\ ideals of shape $(n,\kappa,r,s)$ will be $\vartheta$-efficient if and only if the map $\vartheta$ {\rm (\ref{eqn:mapSigma})} is surjective, which is likely to be the case for general $I$ provided that
\[
     \dim_{\gf}(\operatorname{domain}(\vartheta))\ \geq\ \dim_{\gf}(\operatorname{codomain}(\vartheta)).
\]
This inequality is therefore our precise statement of condition 2.   
\end{description}

\begin{rem} \label{rem:PlausNotPerfect}
    As noted in Remark \ref{rem:genButNotEff}, sufficiently general \SLI\ ideals $I$ of shape $(6,3,2,3)$ are generic and efficient, but not $\vartheta$-efficient.  Indeed, as shown in the associated \emph{Mathematica} notebook, the domain and co-domain of $\vartheta$ have dimensions 90 and 91, respectively, so condition 2 fails in this case, implying that $(6,3,2,3)$ is not a plausible (as defined) shape. This shows that the plausibility criterion is a blunt instrument, incapable of detecting all shapes that support generic \SLI\ ideals.
\end{rem}

\subsection{Implementation and examples} \label{ssec:ImpAndExamplesOfPlausCrit}

Given the preparations in Section \ref{sec:idealsOfShapeNkRs}, the plausibility criterion is straightforward to program; an implementation titled \textbf{genericityIsPlausible} is provided in the notebook \emph{utility functions.nb} mentioned at the start of Section \ref{sec:examples}. Equation (\ref{eqn:PlausibleShapes}) in the introduction lists several plausible shapes (see the notebook \emph{plausible shapes.nb} for the details). 

\subsection{Final observations and a conjecture} \label{ssec:PlausAsympBehavior}

We first explore the second condition of the plausibility criterion more closely.  By (\ref{eqn:QMonDef}) and (\ref{eqn:mapSigma}), we have that 
\[
  \begin{array}{c}
    \dim(\operatorname{domain}(\vartheta)) = n\cdot \lambda \text{ and } \vspace{.05in}\\
    \dim(\operatorname{codomain}(\vartheta)) = |\partial \leadmon\, \cup\, \partial \trailmon| = |\tarmons_{r+1}\, \cup\, \partial \mathcal{O}_{s+1}|;
  \end{array}
\]
therefore, (\ref{eqn:lambdaAndTau}), (\ref{eqn:listOfBdryMons}), and (\ref{eqn:tarmonsRplus1}) yield that condition 2 can be written as follows (recall that we are assuming $2 \leq r < s$):

\[
  \begin{array}{c} 
   n\cdot \left( {n-1+r \choose r} - {\kappa - 1 + r \choose r} \right)\  \geq  {}\vspace{.05in}\\
           {n + r \choose r+1} - {k + r \choose r+1}  + (n-\kappa)\cdot {\kappa - 1 + s \choose s} + { \kappa + s \choose s+1}.
  \end{array}
\]

     We note the following regarding the asymptotic behavior of this inequality when various of the parameters are held constant:
\begin{description}
  \item[Hold $\kappa,r,s$ constant] Since the LHS of the inequalities has dominant term $\frac{n^{r+1}}{r!}$ and the RHS has dominant term $\frac{n^{r+1}}{(r+1)!}$, we see that the inequality holds for all $n>>0$.
  \item[Hold $n,\kappa,r$ constant]  As $s$ increases, we see that the LHS of the inequality is constant and the RHS is increasing, so the inequality will fail for all $s >> 0$.
  \item[Hold $n,r,s$ constant] As $\kappa$ increases (bounded above by $n$, of course), the LHS decreases to $0$ while the RHS is bounded below by ${\kappa + s  \choose s+1 }$ , so there exists $\kappa_0$ $\leq$ $n$ such that the inequality fails for all $\kappa$ $\geq$ $\kappa_0$. 
\end{description}

Next we look more closely at the first condition of the plausibility criterion.  If we hold $\kappa$, $r$, and $s$ constant and let $n$ vary, Lemma \ref{lem:ACountLem} shows that $A_s$ is the most rapidly growing of the sets $A_d$, with dominant term 
\[
    |A_s|\ \approx\ \nu \cdot \frac{n^{r-1}}{(r-1)!}\ \approx\ \frac{n^r}{r!} \cdot \frac{n^{r-1}}{(r-1)!}\ =\ \frac{n^{2r-1}}{r!\cdot (r-1)!}
\]
On the other hand, the number of tangent space relations of degree $s$ is 
\[
    \psi \cdot \tau\ \approx\ \tau \cdot \frac{r}{(r+1)!}\cdot n^{r+1}, \text{ where } \tau = {\kappa-1+s  \choose s } \text{ is independent of }n.
\]

In case $r = 2$, the number of tangent space relations of degree $s$ and $|A_s|$ both grow at the same rate $O(n^{3})$, and the dominant term for the former has the larger coefficient $\tau \cdot \frac{r}{(r+1)!}$ (recall we are assuming $s>r$ and $\kappa$ $>$ $1$, so $\tau$ $\geq$ $s+1$ $\geq$ $4$).  From this it follows that the first condition will be satisfied in degree $s$ (the degree for which satisfaction of condition $1$ is most difficult) for all $n >> 0$.  Hence it is likely that the shape $(n,\kappa,2,s)$ will satisfy condition 1 of the plausibility criterion (as well as condition 2 as seen above) for all $n >> 0$.  This is what leads us to offer 
\medskip

\noindent \textbf{Conjecture \ref{conj:FirstOfTwo}:} Given $r=2$, $s>2$, and $\kappa \geq 2$, the shape $(n,\kappa,2,s)$  is plausible for all $n>>0$.
\medskip

The analogous conjecture for $r>2$ cannot hold: Indeed, if $r$ $>$ $2$, $s > r$, and $\kappa \geq 2$ are fixed, the growth rate $O(n^{2r-1})$ of $|A_s|$ exceeds the growth rate $O(n^{r+1})$, so, as $n$ increases, eventually $|A_s|$ will greatly exceed the number of tangent space relations of degree $s$, thereby falsifying condition 1. Moreover, it appears that for certain choices of $\kappa$, $r$, and $s$, none of the shapes $(n,\kappa,r,s)$ will be plausible; for example, $(n,3,10,11)$ is not plausible for $4$ $\leq$ $n$ $\leq$ $50000$ (at least).  In concluding this paper, we invite the reader to seek further conjectures (or theorems!) regarding families of plausible shapes.

\bibliographystyle{amsplain}
\bibliography{refs}

\end{document}